# Frobenius–Witt cotangent complex for derived rings*

BY ZHOUHANG MAO

*September 7, 2026*

**Abstract**

Shimada recently showed that Frobenius–Witt cotangent complex, the animation of Saito's Frobenius–Witt differentials, vanishes on perfectoid rings, thus it serves as a candidate of "absolute" cotangent complex.

In this article, we give a pullback description of their arithmetic extension, and as a consequence, we give a direct description of Frobenius–Witt cotangent complex, which leads to generalizations to derived rings and animated pre-log rings. This description allows us to compute Frobenius–Witt cotangent complex of derived $\delta$-rings as well.

We also propose a version of Frobenius–Witt cotangent complex relative to derived $\delta$-rings, and establish a vanishing result for prisms, which generalizes vanishing of Frobenius–Witt cotangent complex of perfectoid rings.

Finally, independently of previous considerations, we give a regularity criterion via Frobenius–Witt cotangent complex for $p$-local Noetherian (not necessarily local) rings without $F$-finiteness. We also record the flatness of Frobenius–Witt cotangent complex of valuation rings, with essential ideas due to ChatGPT-6 Astra.

## 1 Introduction

The ring $\mathbb{Z}$ of integers is the initial object in the category of rings, thus it could be understood as an absolute base in algebraic geometry. However, we remark that the affine scheme $\mathrm{Spec}(\mathbb{Z})$ is of dimension 1, not only in the sense of Krull dimension, but also cohomologically. For our purposes, we know that the cotangent complex $L_{\mathbb{F}_p/\mathbb{Z}}$ does not vanish, and for a commutative $\mathbb{F}_p$-algebra $R$, say, the non-vanishing of $L_{\mathbb{F}_p/\mathbb{Z}}$ leads to junk data in the cotangent complex $L_{R/\mathbb{Z}}$, thus we would expect an "absolute" version of cotangent complex, which vanishes on perfect fields of char $p$. Kanau SHIMADA showed in [Shi26] that, the animation of Saito's Frobenius–Witt differentials, *Frobenius–Witt cotangent complex*, serves as a good candidate. More precisely, Frobenius–Witt cotangent complex vanishes on perfectoid rings.

The current article is devoted to a direct description of Frobenius–Witt cotangent complex, which leads to generalizations to (nonconnective) derived rings and animated pre-log rings. The starting point is a pullback description of arithmetic extension in [Shi26].

**Theorem 1.1. (Definition 2.9 and Remarks 2.10 and 2.11)** *Let $R$ be a commutative ring with $\overline{R} := R \otimes^{\mathbb{L}}_{\mathbb{Z}} \mathbb{F}_p$, $M$ an $R/p = \pi_0(\overline{R})$-module, and $m \in M$ an element, viewed as a map $\overline{R} \to M$ of $\overline{R}$-modules. Then the arithmetic extension $W_2(R, M, m)$ fits into a pullback diagram*

$$\begin{array}{ccc} W_2(R, M, m) & \longrightarrow & R \\ \downarrow & & \downarrow \\ \overline{R} \oplus \mathrm{cofib}(m) & \longrightarrow & \overline{R} \oplus \overline{R}[1] \end{array}$$

*of derived rings, where the right vertical map is the composite map*

$$R \xrightarrow{\mathrm{id}_R \otimes 1 \otimes 1} R \otimes^{\mathbb{L}}_{\mathbb{Z}} \mathbb{F}_p \otimes^{\mathbb{L}}_{\mathbb{Z}} \mathbb{F}_p \xrightarrow{\mathrm{Frob}} R \otimes^{\mathbb{L}}_{\mathbb{Z}} \mathbb{F}_p \otimes^{\mathbb{L}}_{\mathbb{Z}} \mathbb{F}_p = \overline{R} \otimes^{\mathbb{L}}_{\mathbb{Z}} \mathbb{F}_p = \overline{R} \oplus \overline{R}[1]$$

*of animated rings, where the map* Frob *is the Frobenius with respect to the last tensor factor $\mathbb{F}_p$.*

This description extends directly to derived rings and animated pre-log rings (Proposition 4.6). By virtue of this description, we get a direct construction of Frobenius–Witt cotangent complex.

*. This article has been written using GNU $\mathrm{T_EX_{MACS}}$ [H+98].

**Theorem 1.2. (Theorems 3.6 and 3.15 and Remark 3.16)** *Let $R$ be a animated ring with $\overline{R} := R \otimes^{\mathbb{L}}_{\mathbb{Z}} \mathbb{F}_p$. The commutative diagram*

$$\begin{array}{ccc} \mathbb{Z} & \longrightarrow & R \\ \downarrow & & \downarrow \\ R & \xrightarrow{\psi_R} & \overline{R} \end{array}$$

*of derived rings, where $\psi_R$ is the composite map $R \xrightarrow{\text{can}} \overline{R} \xrightarrow{\text{Frob}} \overline{R}$, and the right vertical map $R \to \overline{R}$ is the canonical map, induces a map*

$$\partial_R : L_{R/\mathbb{Z}} \otimes^{\mathbb{L}}_{R,\psi_R} \overline{R} \longrightarrow L_{\overline{R}/R} = \overline{R}[1]$$

*of $\overline{R}$-modules, which gives rise to a fiber sequence*

$$\overline{R} \xrightarrow{m_R} \mathrm{fib}(\partial_R) \longrightarrow L_{R/\mathbb{Z}} \otimes^{\mathbb{L}}_{R,\psi_R} \overline{R}.$$

*Then the Frobenius–Witt cotangent complex of $R$ is given by the pointed $\overline{R}$-module $(\mathrm{fib}(\partial_R), m_R)$.*

Again, this description adapts to derived rings and animated pre-log rings (Theorems 4.8 and 4.9). As a consequence of this description, we get a generalization of [Shi26, Prop 6.1] from perfect $\delta$-rings to derived $\delta$-rings.

**Theorem 1.3. (Corollary 3.10 and Remark 3.17)** *Let $R$ be a derived $\delta$-ring with $\overline{R} := R \otimes^{\mathbb{L}}_{\mathbb{Z}} \mathbb{F}_p$. Then its Frobenius–Witt cotangent complex is $\overline{R} \oplus L_{R/\mathbb{Z}} \otimes^{\mathbb{L}}_{R,\psi_R} \overline{R}$.*

Moreover, inspired by this description and [AKN23], we introduce Frobenius–Witt cotangent complex relative to derived $\delta$-rings, which coincides with the absolute one when the base is a perfect $\delta$-ring (Lemma 3.25). We thus get a vanishing of relative Frobenius–Witt cotangent complex for (derived) prisms, generalizing the vanishing of Frobenius–Witt cotangent complex of perfectoid rings.

**Theorem 1.4. (Theorem 3.26)** *Let $(A, I)$ be a (derived) prism, and $R := A/I$ with $\overline{R} := R \otimes^{\mathbb{L}}_{\mathbb{Z}} \mathbb{F}_p$. Then the Frobenius–Witt cotangent complex $FL_{R/A}$ relative to $A$ vanishes.*

Independently of previous developments, we give a regularity criterion via Frobenius–Witt cotangent complex. Unlike relevant results in [Sai22, HJ24, Tak26, Shi26], we assume neither $F$-finiteness nor locality.

**Theorem 1.5. (Theorem 5.2)** *Let $R$ be a $p$-local*[1.1] *Noetherian ring with $\overline{R} := R \otimes^{\mathbb{L}}_{\mathbb{Z}} \mathbb{F}_p$. Then the ring $R$ is regular if and only if its Frobenius–Witt cotangent complex $FL_R$, as an $\overline{R}$-module, is flat.*

Therefore flatness of Frobenius–Witt cotangent complex might be a useful (partial) replacement in the non-Noetherian case, which is further supported by the following theorem, whose proof was developed in discussions with ChatGPT-6 Astra.

**Theorem 1.6. (Theorem A.2)** *Let $V$ be a valuation ring with $\bar{V} := V \otimes^{\mathbb{L}}_{\mathbb{Z}} \mathbb{F}_p$. Then the Frobenius–Witt cotangent complex $FL_V$ of $V$, as a $\bar{V}$-module, is flat.*

**Warning 1.7.** As in [Shi26], Frobenius–Witt cotangent complex is more comparable to Frobenius–Witt differentials “modulo $p$” (cf. [Tak26, Def 3.5]) instead of them on the nose, but if we work with derived $\mathbb{Z}_{(p)}$-algebras, then there is no difference by [Sai22, Lem 1.2(3)].

We work with (nonconnective) derived (commutative) rings, following the theory of Akhil Mathew, based on ideas of Lukas Brantner, as presented in [Rak26, §4]. We denote by DAlg the $\infty$-category of derived rings, and connective derived rings are animated (commutative) rings. Derived $\delta$-rings are understood in the sense of [Hol23, §2.3–§2.4].

1.1. We say that an (animated) ring $R$ is *$p$-local* if $p$ lies in the Jacobson radical $\mathrm{Rad}(\pi_0(R))$.

**Acknowledgments.** The author acknowledges support from the Deutsche Forschungsgemeinschaft (DFG, German Research Foundation) through the Collaborative Research Centre TRR 326 *Geometry and Arithmetic of Uniformized Structures*, project number 444845124. Appendix A is produced under discussions with ChatGPT-6 Astra. We also used ChatGPT-6 Astra for proofreading.

## 2 Arithmetic extension as a pullback

In this section, we give a pull-back description of the arithmetic extension in [Shi26, Def 2.4] (Remark 2.10), which allows us to extend it to derived rings with pointed modules. Technically, we first establish the same result for commutative $\mathbb{F}_p$-algebras (Proposition 2.6), and then extend to arbitrary commutative rings.

Recall that the *arithmetic extension* $W_2(A, M, m)$ of $A$ by $(M, m)$, for a commutative ring $A$, an $A/p$-module $M$ (concentrated in degree 0), and an element $m \in M$, is defined in [Shi26, Def 2.4] to be a ring, with underlying set $A \times M$, and the addition and multiplication is given by

$$\begin{aligned}(a, x) + (b, y) &= \left(a + b, x + y - \sum_{j=1}^{p-1} \frac{1}{p}\binom{p}{j} a^{p-j} b^j m\right),\\ (a, x)(b, y) &= (a b, a^p y + b^p x),\end{aligned}$$

with multiplicative identity given by $(1, 0)$. We note that $(a, 0) + (0, x) = (a, x)$. We start with the special case that $A$ is a commutative $\mathbb{F}_p$-algebra.

**Example 2.1.** Let $A$ be a commutative $\mathbb{F}_p$-algebra. Then by definition, the ring $W_2(A, A, 1)$ coincides with the ring $W_2(A)$ of 2-truncated ($p$-typical) Witt vectors. In this case, by [BS22, Rem 2.5], we know that $W_2(A)$ fits into a Cartesian diagram

$$\begin{array}{ccc} W_2(A) & \longrightarrow & A \\ \downarrow & & \downarrow \\ A & \xrightarrow{\text{can}} & A \otimes_{\mathbb{Z}}^{\mathbb{L}} \mathbb{F}_p \end{array}$$

of derived rings, where the right vertical map is the composite map

$$\psi_A : A \xrightarrow{\text{can}} A \otimes_{\mathbb{Z}}^{\mathbb{L}} \mathbb{F}_p \xrightarrow{\varphi_{A \otimes_{\mathbb{Z}}^{\mathbb{L}} \mathbb{F}_p}} A \otimes_{\mathbb{Z}}^{\mathbb{L}} \mathbb{F}_p.$$

Note that, for any $\mathbb{F}_p$-algebra $A$, its derived modulo $p$ admits a structure of trivial square zero extension.

**Construction 2.2.** Let $A$ be a derived $\mathbb{F}_p$-algebra. Then the derived ring $A \otimes_{\mathbb{Z}}^{\mathbb{L}} \mathbb{F}_p$ can be canonically identified with the square-zero extension $A \oplus A[1]$. Indeed, as

$$A \otimes_{\mathbb{Z}}^{\mathbb{L}} \mathbb{F}_p = A \otimes_{\mathbb{F}_p, 1 \otimes \mathrm{id}_{\mathbb{F}_p}}^{\mathbb{L}} (\mathbb{F}_p \otimes_{\mathbb{Z}}^{\mathbb{L}} \mathbb{F}_p),$$

it suffices to identify the animated $\mathbb{F}_p$-algebra $\mathbb{F}_p \otimes_{\mathbb{Z}}^{\mathbb{L}} \mathbb{F}_p$, where the $\mathbb{F}_p$-algebra structure is given by the second factor, as a square-zero extension $\mathbb{F}_p \oplus \mathbb{F}_p[1]$. We identify the first factor $\mathbb{F}_p$ as the pushout $\mathbb{Z} \otimes_{0 \leftarrow t, \mathbb{Z}[t], t \mapsto p}^{\mathbb{L}} \mathbb{Z}$ of animated rings, and thus the animated $\mathbb{F}_p$-algebra $\mathbb{F}_p \otimes_{\mathbb{Z}}^{\mathbb{L}} \mathbb{F}_p$ can be identified with $\mathbb{F}_p \otimes_{\mathbb{F}_p[t]}^{\mathbb{L}} \mathbb{F}_p = \mathrm{LSym}_{\mathbb{F}_p}(\mathbb{F}_p[1])$, which gives rise to a square-zero extension structure.

Example 2.1 hints the following definition.

**Definition 2.3.** *Let $A$ be a derived $\mathbb{F}_p$-algebra, and $M$ an $A$-module spectrum with a map $m: A \to M$ of $A$-modules. The map $m$ gives rise to a map $\mathrm{cofib}(m) \to A[1]$ of $A$-modules. Then we define $\overline{W}_2(A, M, m)$ to be the pullback of the diagram*

$$\begin{array}{ccc} & & A \\ & & \downarrow{\scriptstyle \psi_A} \\ A \oplus \mathrm{cofib}(m) & \longrightarrow & A \oplus A[1] \simeq A \otimes_{\mathbb{Z}}^{\mathbb{L}} \mathbb{F}_p \end{array}$$

*of derived rings. As the map $\psi_A$ is not a map of derived $\mathbb{F}_p$-algebras in general, this diagram is **not** a diagram of derived $\mathbb{F}_p$-algebras. We denote by* res *the canonical map $\overline{W}_2(A, M, m) \to A$.*

**Remark 2.4.** Let $A$ be a derived $\mathbb{F}_p$-algebra, and $M$ an $A$-module spectrum with a map $m: A \to M$ of $A$-modules. By definition, we have a Cartesian diagram

$$\begin{array}{ccc} \overline{W}_2(A, M, m) & \longrightarrow & A \\ \downarrow & & \downarrow \\ A \oplus \mathrm{cofib}(m) & \longrightarrow & A \otimes_{\mathbb{Z}}^{\mathbb{L}} \mathbb{F}_p \end{array}$$

of $\mathbb{Z}$-module spectra, which induces an equivalence of horizontal fibers, and thus we get a fiber sequence

$$M \longrightarrow \overline{W}_2(A, M, m) \xrightarrow{\mathrm{res}} A$$

in $D(\mathbb{Z})$, which is functorial in $(A, M, m)$. We denote by $V$ the first map $M \to \overline{W}_2(A, M, m)$. In particular, if both $A$ and $M$ are concentrated in degree 0, then so is the $\mathbb{Z}$-module spectrum $\overline{W}_2(A, M, m)$, which implies that the derived ring $\overline{W}_2(A, M, m)$ is a classical commutative ring.

**Construction 2.5.** Let $A$ be a derived $\mathbb{F}_p$-algebra, and $M$ an $A$-module spectrum with a map $m: A \to M$ of $A$-modules. Apply Remark 2.4 to the map $(A, A, 1) \to (A, M, m)$ induced by the map $m$, we get a map

$$\begin{array}{ccccc} A & \xrightarrow{V} & W_2(A) & \xrightarrow{\mathrm{res}} & A \\ \downarrow & & \downarrow & & \| \\ M & \xrightarrow{V} & \overline{W}_2(A, M, m) & \xrightarrow{\mathrm{res}} & A \end{array} \tag{2.1}$$

of fiber sequences in $D(\mathbb{Z})$. Since the rightmost vertical map is an equivalence, the left square is Cartesian in $D(\mathbb{Z})$. When both $A$ and $M$ are connective, we get a multiplicative section $\Omega^\infty A \to \Omega^\infty \overline{W}_2(A, M, m)$ of the canonical map res$: \overline{W}_2(A, M, m) \to A$ given by the composite map

$$[-]_{(A,M,m)}: \Omega^\infty A \xrightarrow{[-]_A} \Omega^\infty W_2(A) \longrightarrow \Omega^\infty \overline{W}_2(A, M, m),$$

where the first map $[-]_A$ is the Teichmüller map, being a multiplicative section of the restriction map res$: W_2(A) \to A$. Consequently, we get an equivalence

$$\Omega^\infty A \times \Omega^\infty M \longrightarrow \Omega^\infty \overline{W}_2(A, M, m)$$

of anima, informally written as $(a, m) \mapsto ([a]_{(A,M,m)} + V(m))$.

**Proposition 2.6.** *Let $A$ be a commutative $\mathbb{F}_p$-algebra, and $M$ an $A$-module concentrated in degree 0, and $m \in M$ an element inducing a map $m: A \to M$ of $A$-modules. Then the isomorphism*

$$W_2(A, M, m) = A \times M \to \overline{W}_2(A, M, m)$$

*of sets (we invoke Remark 2.4 to see that they are sets, so we omit notations $\Omega^\infty$) is compatible with the ring structures on both sides.*

**Proof.** As $(1, 0)$ is the unit of $W_2(A)$, this map preserves multiplicative identities. The additivity is verified by the left Cartesian square of (2.1) in $D(\mathbb{Z})$ in Construction 2.5.

Now we verify the multiplicativity. Note that the map $A \oplus \mathrm{cofib}(m) \to A \oplus A[1]$ of derived rings in Definition 2.3 can be realized as a square-zero extension, thus so can the (surjective) map res$: \overline{W}_2(A, M, m) \to A$ of rings. Given the additivity, the multiplicativity is completely determined by the $A$-module structure on $M$, which is given by forgetting the $A \otimes_{\mathbb{Z}}^{\mathbb{L}} \mathbb{F}_p$-module structure on $M$ along the map $\psi_A$. Now $M$ is concentrated in degree 0, so it is completely the $\pi_0(A \otimes_{\mathbb{Z}}^{\mathbb{L}} \mathbb{F}_p) = A$-module structure on $M$, and the composite $A \xrightarrow{\psi_A} A \otimes_{\mathbb{Z}}^{\mathbb{L}} \mathbb{F}_p \xrightarrow{\tau_{\leq 0}} \pi_0(A \otimes_{\mathbb{Z}}^{\mathbb{L}} \mathbb{F}_p) = A$ is simply the Frobenius map $\varphi_A$. The multiplicativity then follows. □

Now we extend to the case that $A$ is a general derived ring (in place of a derived $\mathbb{F}_p$-algebra).

**Remark 2.7.** Let $A$ be a commutative ring, and $M$ an $A/p$-module concentrated in degree 0, and $m \in M$ an element. Then we have a map

$$\begin{array}{ccccccccc} 0 & \longrightarrow & M & \longrightarrow & W_2(A, M, m) & \longrightarrow & A & \longrightarrow & 0 \\ & & \| & & \downarrow & & \downarrow & & \\ 0 & \longrightarrow & M & \longrightarrow & W_2(A/p, M, m) & \longrightarrow & A/p & \longrightarrow & 0 \end{array}$$

of short exact sequences, where the rightmost vertical map is the canonical map $A \to A/p$. It follows that the right square is Cartesian in $D(\mathbb{Z})$ (not just in abelian groups).

We need the following simple observation to see that the difference between $\pi_0(A \otimes_{\mathbb{Z}}^{\mathbb{L}} \mathbb{F}_p) = A/p$ and $A \otimes_{\mathbb{Z}}^{\mathbb{L}} \mathbb{F}_p$ is insignificant.

**Remark 2.8.** Let $f : A \to B$ be a map of derived $\mathbb{F}_p$-algebras, and $M$ a $B$-module, and $m : B \to M$ a map of $B$-modules. Then by Remark 2.4, the square

$$\begin{array}{ccc} \overline{W}_2(A, M, m \circ f) & \longrightarrow & A \\ \downarrow & & \downarrow \\ \overline{W}_2(B, M, m) & \longrightarrow & B \end{array}$$

of derived rings is Cartesian.

We come to the following definition.

**Definition 2.9.** *Let $A$ be a derived ring with $\overline{A} := A \otimes_{\mathbb{Z}}^{\mathbb{L}} \mathbb{F}_p$, and $M$ an $\overline{A}$-module, and $m : \overline{A} \to M$ a map of $\overline{A}$-modules. Then we define the* derived arithmetic extension $W_2^{\mathrm{der}}(A, M, m)$ *(which will be denoted by $W_2(A, M, m)$ after Remark 2.11) of $A$ by $M$ to be the pullback of the diagram*

$$\begin{array}{ccc} & & A \\ & & \downarrow \\ \overline{W}_2(\overline{A}, M, m) & \longrightarrow & A \otimes_{\mathbb{Z}}^{\mathbb{L}} \mathbb{F}_p \end{array}$$

*of derived rings, where the right vertical map is the canonical map* $\mathrm{can} : A \to \overline{A}$ *(not the Frobenius). We denote by* res *the canonical map $W_2^{\mathrm{der}}(A, M, m) \to A$.*

**Remark 2.10.** Let $A$ be a derived ring with $\overline{A} := A \otimes_{\mathbb{Z}}^{\mathbb{L}} \mathbb{F}_p$, and $M$ an $\overline{A}$-module, and $m : \overline{A} \to M$ a map of $\overline{A}$-modules. Unraveling definitions, we see that the derived ring $W_2^{\mathrm{der}}(A, M, m)$ fits canonically into a Cartesian diagram

$$\begin{array}{ccc} W_2^{\mathrm{der}}(A, M, m) & \longrightarrow & A \\ \downarrow & & \downarrow \\ \overline{A} \oplus \mathrm{cofib}(m) & \longrightarrow & \overline{A} \oplus \overline{A}[1] \end{array},$$

where the right vertical map is the composite $A \xrightarrow{\mathrm{can}} \overline{A} \xrightarrow{\psi_{\overline{A}}} \overline{A} \otimes_{\mathbb{Z}}^{\mathbb{L}} \mathbb{F}_p = \overline{A} \oplus \overline{A}[1]$, and the bottom horizontal map is induced by the boundary map $\mathrm{cofib}(m) \to \overline{A}[1]$.

**Remark 2.11.** Let $A$ be a commutative ring, and $M$ an $A/p$-module concentrated in degree 0, and $m \in M$ an element. Remark 2.8 shows that the square

$$\begin{array}{ccc} \overline{W}_2(A \otimes_{\mathbb{Z}}^{\mathbb{L}} \mathbb{F}_p, M, m) & \longrightarrow & A \otimes_{\mathbb{Z}}^{\mathbb{L}} \mathbb{F}_p \\ \downarrow & & \downarrow \\ \overline{W}_2(A/p, M, m) & \longrightarrow & A/p \end{array}$$

of derived rings is Cartesian. Thus by Remark 2.7, it follows from Proposition 2.6 that the functor $W_2^{\mathrm{der}}$ defined in Definition 2.9 coincides with the functor $W_2$ in [Shi26, Def 2.4] reviewed above when derived rings and modules in question are concentrated in degree 0, thus we will denote $W_2^{\mathrm{der}}$ by $W_2$ in the sequel.

**Remark 2.12.** Let $A$ be a derived $\mathbb{F}_p$-algebra, $M$ an $A$-module, and $m: A \to M$ a map of $A$-modules. The multiplication map $\mu_A: A \otimes_{\mathbb{Z}}^{\mathbb{L}} \mathbb{F}_p \to A$ endows $M$ an $A \otimes_{\mathbb{Z}}^{\mathbb{L}} \mathbb{F}_p$-module structure. Now we examine the diagram

$$\begin{array}{ccc} W_2(A, M, m) & \longrightarrow & A \\ \downarrow & & \downarrow \\ \overline{W}_2(A \otimes_{\mathbb{Z}}^{\mathbb{L}} \mathbb{F}_p, M, m \circ \mu_A) & \longrightarrow & A \otimes_{\mathbb{Z}}^{\mathbb{L}} \mathbb{F}_p \\ \downarrow & & \downarrow \\ \overline{W}_2(A, M, m) & \longrightarrow & A \end{array}$$

of derived rings. By definition, the top square is Cartesian. By Remark 2.8, the bottom square is Cartesian as well, thus so is the outer square. Since the composite map

$$A \xrightarrow{\text{can}} A \otimes_{\mathbb{Z}}^{\mathbb{L}} \mathbb{F}_p \xrightarrow{\mu_A} A$$

is an equivalence, it follows that the map $W_2(A, M, m) \to \overline{W}_2(A, M, m)$ is an equivalence as well.

**Remark 2.13.** Let $A$ be a derived ring with $\overline{A} := A \otimes_{\mathbb{Z}}^{\mathbb{L}} \mathbb{F}_p$, and $M$ an $\overline{A}$-module, and $m: \overline{A} \to M$ a map of $\overline{A}$-modules. By Definition 2.9 and Remark 2.4, we have a fiber sequence

$$M \longrightarrow W_2(A, M, m) \xrightarrow{\text{res}} A$$

in $D(\mathbb{Z})$. We denote by $V$ the first map $M \to W_2(A, M, m)$. This shows that the functor $W_2$ preserves sifted colimits, and in particular, if both $A$ and $M$ are connective, then so is $W_2(A, M, m)$, and this coincides canonically with $W_2^{\text{an}}(A, M, m)$ in [Shi26, Def 4.1].

**Remark 2.14.** Let $f: A \to B$ be a map of derived rings with $\overline{A} := A \otimes_{\mathbb{Z}}^{\mathbb{L}} \mathbb{F}_p$, $\overline{B} := B \otimes_{\mathbb{Z}}^{\mathbb{L}} \mathbb{F}_p$, and $\overline{f} := f \otimes_{\mathbb{Z}}^{\mathbb{L}} \mathbb{F}_p$, and $M$ a $\overline{B}$-module, and $m: \overline{B} \to M$ a map of $B$-modules. Then by Remark 2.4, the square

$$\begin{array}{ccc} W_2(A, M, m \circ \overline{f}) & \longrightarrow & A \\ \downarrow & & \downarrow \\ W_2(B, M, m) & \longrightarrow & B \end{array}$$

of derived rings is Cartesian.

## 3 Frobenius–Witt cotangent complex relative to $\delta$-rings

In this section, thanks to the pullback description in Remark 2.10, we give a rather direct description of Frobenius–Witt cotangent complex (Theorems 3.6 and 3.15), which extends directly to derived rings. As a consequence, we deduce a computation of Frobenius–Witt cotangent complex of derived $\delta$-rings (Corollary 3.10 and Remark 3.17), generalizing the result [Shi26, Prop 6.1] for perfect $\delta$-rings. Inspired by this description, we also introduce the concept of Frobenius–Witt cotangent complex relative to derived $\delta$-rings (Definition 3.19), and deduce a generalization of vanishing of Frobenius–Witt cotangent complex of perfectoid rings (Theorem 3.26).

**Construction 3.1.** Let $R$ be a derived $\mathbb{F}_p$-algebra. The map $\psi_R: R \to R \otimes_{\mathbb{Z}}^{\mathbb{L}} \mathbb{F}_p = R \oplus R[1]$ of derived rings (here we are invoking Construction 2.2) is then classified by a map $\overline{\partial}_R: L_{R/\mathbb{Z}} \otimes_{R, \varphi_R}^{\mathbb{L}} R \to R[1]$ of $R$-modules, where we use the fact that the composite map

$$R \xrightarrow{\psi_R} R \otimes_{\mathbb{Z}}^{\mathbb{L}} \mathbb{F}_p \xrightarrow{\mu_R} R$$

coincides with the Frobenius map $\varphi_R$. This gives rise to a map $\overline{m}_R: R \to \operatorname{fib}(\overline{\partial}_R)$ of $R$-modules, along with a fiber sequence

$$R \xrightarrow{\overline{m}_R} \operatorname{fib}(\overline{\partial}_R) \longrightarrow L_{R/\mathbb{Z}} \otimes_{R, \varphi_R}^{\mathbb{L}} R$$

of $R$-modules.

**Lemma 3.2.** *The map $\overline{\partial}_{\mathbb{F}_p}$ of $\mathbb{F}_p$-modules in Construction 3.1 is an equivalence.*

**Proof.** Note that the source and the target of $\overline{\partial}_{\mathbb{F}_p}$ are both $\mathbb{F}_p[1]$, thus $\overline{\partial}_{\mathbb{F}_p}$ is either 0 or an equivalence. It suffices to exclude the possibility of $\overline{\partial}_{\mathbb{F}_p}=0$. Indeed, if it were zero, then the diagram $\mathbb{F}_p \xrightarrow{\text{can}} \mathbb{F}_p \oplus \mathbb{F}_p[1] \xleftarrow{\psi_{\mathbb{F}_p}} \mathbb{F}_p$ would be the trivial square-zero extension $\mathbb{F}_p \oplus \mathbb{F}_p$, while it is actually $W_2(\mathbb{F}_p) = \mathbb{Z}/p^2$. □

**Construction 3.3.** Let $R$ be a derived ring. Apply Construction 3.1 to the derived $\mathbb{F}_p$-algebra $\overline{R} = R \otimes^{\mathbb{L}}_{\mathbb{Z}} \mathbb{F}_p$, we get a composite map

$$\partial_R : L_{R/\mathbb{Z}} \otimes^{\mathbb{L}}_{R,\psi_R} \overline{R} \longrightarrow L_{\overline{R}/\mathbb{Z}} \otimes^{\mathbb{L}}_{\overline{R},\varphi_{\overline{R}}} \overline{R} \xrightarrow{\overline{\partial}_{\overline{R}}} \overline{R}[1]$$

of $\overline{R}$-modules classifying the composite map $R \xrightarrow{\text{can}} \overline{R} \xrightarrow{\psi_R} \overline{R} \otimes^{\mathbb{L}}_{\mathbb{Z}} \mathbb{F}_p = \overline{R} \oplus \overline{R}[1]$ of derived rings, which gives rise to a map $m_R : \overline{R} \to \mathrm{fib}(\partial_R)$ of $\overline{R}$-modules, along with a fiber sequence

$$\overline{R} \xrightarrow{m_R} \mathrm{fib}(\partial_R) \longrightarrow L_{R/\mathbb{Z}} \otimes^{\mathbb{L}}_{R,\psi_R} \overline{R}$$

of $\overline{R}$-modules.

**Lemma 3.4.** *The map $\partial_{\mathbb{F}_p}$ in Construction 3.3 is an equivalence.*

**Proof.** Let $R := \mathbb{F}_p$, and $\overline{R} := R \otimes^{\mathbb{L}}_{\mathbb{Z}} \mathbb{F}_p$. The map $\overline{R} \to R$ of derived $\mathbb{F}_p$-algebras gives rise to a commutative diagram

$$\begin{array}{ccccccc} R & \xrightarrow{\text{can}} & \overline{R} & \xrightarrow{\psi_{\overline{R}}} & \overline{R} \otimes^{\mathbb{L}}_{\mathbb{Z}} \mathbb{F}_p & = & \overline{R} \oplus \overline{R}[1] \\ & & \downarrow & & \downarrow & & \downarrow \\ & & R & \xrightarrow{\psi_R} & R \otimes^{\mathbb{L}}_{\mathbb{Z}} \mathbb{F}_p & = & R \oplus R[1] \end{array}$$

of derived rings (the commutativity of the rightmost square uses the map $\overline{R} \to R$ of derived $\mathbb{F}_p$-algebras, not merely as a map of derived rings). It follows that $\overline{\partial}_R = \partial_R \otimes^{\mathbb{L}}_{\overline{R}} R$. The result then follows from Lemma 3.2 along with the conservativity of the functor $(-) \otimes^{\mathbb{L}}_{\overline{R}} R : D(\overline{R})^{\omega} \to D(R)^{\omega}$ on perfect complexes. □

As in [Shi26, Def 4.5], we define Frobenius–Witt derivations in the derived context.

**Definition 3.5.** *Let $R$ be a derived ring with $\overline{R} := R \otimes^{\mathbb{L}}_{\mathbb{Z}} \mathbb{F}_p$, and $M$ an $\overline{R}$-module, and $m : \overline{R} \to M$ a map of $\overline{R}$-modules. Then the anima* $\mathrm{FWDer}(R,(M,m))$ *of* (pointed) Frobenius–Witt derivations *is defined to be*

$$\mathrm{Hom}_{\mathrm{DAlg}_{/R}}(R, W_2(R,M,m)),$$

*where the map $W_2(R,M,m) \to R$ is the canonical map* res *in Definition 2.9.*

The following statement shows that $(\mathrm{fib}(\partial_R), m_R)$ coincides with Frobenius–Witt cotangent complex in [Shi26, Prop 4.8] for animated $\mathbb{Z}_{(p)}$-algebras $R$, and in particular, it could be viewed as the definition of Frobenius–Witt cotangent complex for derived $\mathbb{Z}_{(p)}$-algebras. In particular, this gives rise to the vanishing of Frobenius–Witt cotangent complex [Shi26, Lem 6.3] without computations.

**Theorem 3.6.** *Let $R$ be a derived ring with $\overline{R} := R \otimes^{\mathbb{L}}_{\mathbb{Z}} \mathbb{F}_p$. Then the functor*

$$\begin{array}{ccc} D(\overline{R})_{\overline{R}/} & \longrightarrow & \mathrm{An} \\ (M,m) & \longmapsto & \mathrm{FWDer}(R,(M,m)) \end{array}$$

*is co-represented by* $(\mathrm{fib}(\partial_R), m_R)$ *in Construction 3.3.*

**Proof.** Let $M$ be an $\overline{R}$-module, and $m : \overline{R} \to M$ a map of $R$-modules. Unraveling definitions, by Remark 2.10, we see that the anima $\mathrm{FWDer}(R,(M,m))$ is given by the anima of a map

$$\alpha : L_{R/\mathbb{Z}} \otimes^{\mathbb{L}}_{R,\psi_R} \overline{R} \longrightarrow \mathrm{cofib}(m)$$

along with a homotopy between the map $\partial_R : L_{R/\mathbb{Z}} \otimes^{\mathbb{L}}_{R,\psi_R} \overline{R} \to \overline{R}[1]$ of $\overline{R}$-modules in Construction 3.3 and the composite map

$$L_{R/\mathbb{Z}} \otimes^{\mathbb{L}}_{R,\psi_R} \overline{R} \xrightarrow{\alpha} \mathrm{cofib}(m) \longrightarrow \overline{R}[1]$$

of $\overline{R}$-modules, i.e. the anima of factorizations of $\partial_R$ through $\mathrm{cofib}(m) \to \overline{R}[1]$ in $D(\overline{R})$. On the other hand, as the $\infty$-category $D(\overline{R})$ is stable, we have an equivalence

$$\begin{array}{ccc} D(\overline{R})_{\overline{R}/} & \longrightarrow & D(\overline{R})_{/\overline{R}[1]} \\ (\beta : \overline{R} \to M) & \longmapsto & (\mathrm{cofib}(\beta) \to \overline{R}[1]) \end{array}$$

of $\infty$-categories, and in particular,

$$\mathrm{Hom}_{D(\overline{R})_{\overline{R}/}}((\mathrm{fib}(\partial_R), m_R), (M, m)) = \mathrm{Hom}_{D(\overline{R})_{/\overline{R}[1]}}(\partial_R, \mathrm{cofib}(m) \to \overline{R}[1])$$

where the map $\mathrm{cofib}(m) \to \overline{R}[1]$ is the boundary map. The result then follows. □

Such a description would be helpful even if we restrict to animated rings or even classical rings, as long as we have a way to control its Frobenius. For example, the Frobenius map of $\delta$-rings are well behaved. We have the following simple observation.

**Lemma 3.7.** *Let $R$ be a derived $\delta$-ring, and $S$ a derived $\mathbb{F}_p$-algebra. Then the Frobenius map $\varphi_{R \otimes^{\mathbb{L}}_{\mathbb{Z}} S} : R \otimes^{\mathbb{L}}_{\mathbb{Z}} S \to R \otimes^{\mathbb{L}}_{\mathbb{Z}} S$ can be canonically identified with the map $\varphi_R \otimes^{\mathbb{L}}_{\mathbb{Z}} \varphi_S$.*

**Proof.** Let $\overline{R}$ be the derived $\mathbb{F}_p$-algebra $R \otimes^{\mathbb{L}}_{\mathbb{Z}} \mathbb{F}_p$. Then

$$\begin{aligned} \varphi_R \otimes^{\mathbb{L}}_{\mathbb{Z}} \varphi_S &= \varphi_{\overline{R}} \otimes^{\mathbb{L}}_{\mathbb{F}_p} \varphi_S \\ &= \varphi_{\overline{R} \otimes^{\mathbb{L}}_{\mathbb{F}_p} S} \\ &= \varphi_{R \otimes^{\mathbb{L}}_{\mathbb{Z}} S}, \end{aligned}$$

where the first equality follows from the fact that $\varphi_R$ is a derived Frobenius lift. □

**Corollary 3.8.** *Let $R$ be a derived $\delta$-ring. Then the map $\partial_R$ is canonically null-homotopic.*

**Proof.** By Lemma 3.7 (with $S = \mathbb{F}_p \otimes^{\mathbb{L}}_{\mathbb{Z}} \mathbb{F}_p$), we see that the composite map

$$R \xrightarrow{\mathrm{can}} \overline{R} \xrightarrow{\psi_R} \overline{R} \otimes^{\mathbb{L}}_{\mathbb{Z}} \mathbb{F}_p = \overline{R} \oplus \overline{R}[1]$$

of derived rings can be rewritten as the composite of the map

$$\alpha : R \xrightarrow{\mathrm{id}_R \otimes 1 \otimes 1} R \otimes^{\mathbb{L}}_{\mathbb{Z}} \mathbb{F}_p \otimes^{\mathbb{L}}_{\mathbb{Z}} \mathbb{F}_p \xrightarrow{\mathrm{id}_R \otimes \varphi_{\mathbb{F}_p \otimes^{\mathbb{L}}_{\mathbb{Z}} \mathbb{F}_p}} \overline{R} \otimes^{\mathbb{L}}_{\mathbb{Z}} \mathbb{F}_p = \overline{R} \oplus \overline{R}[1]$$

of derived rings with the map

$$\varphi_{\overline{R}} \otimes^{\mathbb{L}}_{\mathbb{Z}} \mathbb{F}_p : \overline{R} \otimes^{\mathbb{L}}_{\mathbb{Z}} \mathbb{F}_p \xrightarrow{(\varphi_R \otimes^{\mathbb{L}}_{\mathbb{Z}} \mathbb{F}_p) \otimes^{\mathbb{L}}_{\mathbb{Z}} \mathbb{F}_p} \overline{R} \otimes^{\mathbb{L}}_{\mathbb{Z}} \mathbb{F}_p$$

of derived $\mathbb{F}_p$-algebras. Note that the map $\beta : L_{R/\mathbb{Z}} \otimes^{\mathbb{L}}_{R,\mathrm{can}} \overline{R} \to \overline{R}[1]$ classifying the map $\alpha$ into the trivial square-zero extension $\overline{R} \oplus \overline{R}[1]$ is canonically null-homotopic, as the initiality of $\mathbb{Z}$ in DAlg implies that any map $\mathbb{Z} \to \mathbb{F}_p \otimes \mathbb{F}_p = \mathbb{F}_p \oplus \mathbb{F}_p[1]$ of derived rings coincides with the composite map $\mathbb{Z} \to \mathbb{F}_p \hookrightarrow \mathbb{F}_p \oplus \mathbb{F}_p[1]$ of derived rings. It follows that $\partial_R = \beta \otimes^{\mathbb{L}}_{R,\varphi_{\overline{R}}} \overline{R}$ is canonically null-homotopic as well. □

**Warning 3.9.** Let $R$ be a derived ring. Then the composite map

$$R \otimes^{\mathbb{L}}_{\mathbb{Z}} \mathbb{F}_p \otimes^{\mathbb{L}}_{\mathbb{Z}} \mathbb{F}_p \xrightarrow{\mathrm{id}_R \otimes \varphi_{\mathbb{F}_p \otimes^{\mathbb{L}}_{\mathbb{Z}} \mathbb{F}_p}} R \otimes^{\mathbb{L}}_{\mathbb{Z}} \mathbb{F}_p \otimes^{\mathbb{L}}_{\mathbb{Z}} \mathbb{F}_p \xrightarrow{\varphi_{R \otimes^{\mathbb{L}}_{\mathbb{Z}} \mathbb{F}_p} \otimes \mathrm{id}_{\mathbb{F}_p}} R \otimes^{\mathbb{L}}_{\mathbb{Z}} \mathbb{F}_p \otimes^{\mathbb{L}}_{\mathbb{Z}} \mathbb{F}_p$$

of derived $\mathbb{F}_p$-algebras is well-defined. In the proof of Corollary 3.8, we show that, a $\delta$-structure on $R$ gives rise to an identification of this composite map with the Frobenius map $\varphi_{R \otimes^{\mathbb{L}}_{\mathbb{Z}} \mathbb{F}_p \otimes^{\mathbb{L}}_{\mathbb{Z}} \mathbb{F}_p}$ via Lemma 3.7. The $\delta$-structure here is crucial. Actually, when $R = \mathbb{F}_p$, there is no such identification, as if there were, then the proof of Corollary 3.8 would imply that $\partial_R = 0$, but we have already seen in Lemma 3.4 that $\partial_{\mathbb{F}_p}$ is an automorphism of a non-zero spectrum.

It follows from Corollary 3.8 that

**Corollary 3.10.** *Let $R$ be a derived $\delta$-ring, and $\overline{R} := R \otimes^{\mathbb{L}}_{\mathbb{Z}} \mathbb{F}_p$. Then the Frobenius–Witt cotangent complex of $R$ is $\overline{R} \oplus L_{R/\mathbb{Z}} \otimes^{\mathbb{L}}_{R,\psi_R} \overline{R}$, with $m$ being the inclusion $\overline{R} \hookrightarrow \overline{R} \oplus L_{R/\mathbb{Z}} \otimes^{\mathbb{L}}_{R,\psi_R} \overline{R}$.*

**Example 3.11.** Let $R$ be a perfect $\delta$-ring. Then $L_{R/\mathbb{Z}} \otimes^{\mathbb{L}}_{R,\mathrm{can}} \overline{R} = L_{\overline{R}/\mathbb{F}_p} = 0$. By Corollary 3.10 (and the fact that $\psi_R = \varphi_{\overline{R}} \circ \mathrm{can}$), the Frobenius–Witt cotangent complex of $R$ is $\overline{R}$ with $m = 1$, recovering [Shi26, Prop 6.1]. Thus Corollary 3.10 is a generalization of it to non-perfect (derived) $\delta$-rings.

We now give a simple description of the Frobenius map $\varphi_{\mathbb{F}_p \otimes^{\mathbb{L}}_{\mathbb{Z}} \mathbb{F}_p}$, which simplifies the map $\partial_R$ in Construction 3.3 a lot. We first explain how to get a filtered Frobenius map on filtered derived $\mathbb{F}_p$-algebras. One may compare this construction with [AMMN22, Ex A.11].

**Construction 3.12. (Filtered Frobenius)** We construct filtered Frobenius map on $(\mathbb{N}, \geq)$-filtered derived $\mathbb{F}_p$-algebras as in [Hol23, Cons 2.4.1]. Indeed, given an $(\mathbb{N}, \geq)$-filtered commutative $\mathbb{F}_p$-algebra $\mathrm{Fil}\, A$, we simply assign $x \in \mathrm{Fil}^m A$ to $x^p \in \mathrm{Fil}^{mp} A$. This gives rise to a map

$$\varphi_{\mathrm{Fil}\, A} : \mathrm{Fil}\, A \longrightarrow (\mathrm{Fil}^{mp} A)_{m \in \mathbb{N}} \tag{3.1}$$

of $(\mathbb{N}, \geq)$-filtered commutative $\mathbb{F}_p$-algebras. We restrict such construction to free filtered commutative $\mathbb{F}_p$-algebras generated on a finitely filtered finite set, and argue as in [Hol23, Cons 2.4.1], we see that this construction is right-extendable. This extends the filtered Frobenius map $\varphi_{\mathrm{Fil}\, A}$ in (3.1) to $(\mathbb{N}, \geq)$-filtered derived $\mathbb{F}_p$-algebras $\mathrm{Fil}\, A$.

Now we observe that, the Frobenius map on a square-zero extension factors through the extended algebra.

**Construction 3.13.** Let $R$ be a derived $\mathbb{F}_p$-algebra, $M$ an $R$-module, and $S$ a square-zero extension of $R$ by $M$ within $D(\mathbb{F}_p)$. By [Mag24, Prop 3.3.17], we may identify $S$ with an $(\mathbb{N}, \geq)$-filtered derived $\mathbb{F}_p$-algebra $\mathrm{Fil}\, S$ concentrated in weights $[0,1]$, such that $\mathrm{Fil}^0 S = S$, $\mathrm{Fil}^1 S = M$, and $\mathrm{gr}^0 S = R$. We now apply Construction 3.12 to the filtered derived $\mathbb{F}_p$-algebra $\mathrm{Fil}\, S$, obtaining the filtered Frobenius map $\varphi_{\mathrm{Fil}\, S}$. Note that the filtered derived $\mathbb{F}_p$-algebra $(\mathrm{Fil}^{mp} S)_{m \in \mathbb{N}}$ is concentrated in weights $\{0\}$. It follows that the Frobenius map $\varphi_S : S \to S$ coincides with the composite

$$\mathrm{Fil}^0 S \longrightarrow \mathrm{gr}^0 S \xrightarrow{\varphi_{\mathrm{Fil}\, S}} \mathrm{gr}^0 (\mathrm{Fil}^{mp} S)_{m \in \mathbb{N}} \xleftarrow{\simeq} \mathrm{Fil}^0 S.$$

This gives rise to a factorization $S \to R \to S$ of the Frobenius map $\varphi_S$. In particular, we get a factorization

$$\mathbb{F}_p \otimes^{\mathbb{L}}_{\mathbb{Z}} \mathbb{F}_p \xrightarrow{\mathrm{mult}} \mathbb{F}_p \xrightarrow{1 \otimes \mathrm{id}} \mathbb{F}_p \otimes^{\mathbb{L}}_{\mathbb{Z}} \mathbb{F}_p$$

of the Frobenius map $\varphi_{\mathbb{F}_p \otimes^{\mathbb{L}}_{\mathbb{Z}} \mathbb{F}_p}$, where the $\mathbb{F}_p$-algebra structure on each tensor product is given by right factor.

**Remark 3.14.** Construction 3.13 gives rise to another proof of Lemma 3.2. Indeed, the map $\bar{\partial}_{\mathbb{F}_p}$ classifies the map $\mathbb{F}_p \xrightarrow{1 \otimes \mathrm{id}_{\mathbb{F}_p}} \mathbb{F}_p \otimes^{\mathbb{L}}_{\mathbb{Z}} \mathbb{F}_p = \mathbb{F}_p \oplus \mathbb{F}_p[1]$ (cf. Construction 2.2). This is the base change of the map

$$\mathbb{Z} \xrightarrow{1 \otimes \mathrm{id}_{\mathbb{Z}}} \mathbb{Z} \otimes^{\mathbb{L}}_{\mathbb{Z}[t]} \mathbb{Z} = \mathbb{Z} \oplus \mathbb{Z}[1]$$

of animated $\mathbb{Z}[t]$-algebras along $\mathbb{Z}[t] \to \mathbb{Z}, t \mapsto p$, and thus it suffices to see that the classifying map $\gamma : L_{\mathbb{Z}/\mathbb{Z}[t]} \to \mathbb{Z}[1]$ is an equivalence. The coCartesian square

$$\begin{array}{ccc} \mathbb{Z}[t] & \longrightarrow & \mathbb{Z} \\ \downarrow & & \downarrow \\ \mathbb{Z} & \longrightarrow & \mathbb{Z} \otimes^{\mathbb{L}}_{\mathbb{Z}[t]} \mathbb{Z} \end{array}$$

in DAlg is LSym of the (co)Cartesian square

$$\begin{array}{ccc} \mathbb{Z} & \xrightarrow{\alpha} & 0 \\ \downarrow & & \downarrow \\ 0 & \longrightarrow & \mathbb{Z}[1] \end{array}$$

in $D(\mathbb{Z})$, and thus the map $L_{\mathbb{Z}/\mathbb{Z}[t]} \to \mathbb{Z}[1]$ is actually induced by the map $\beta: \mathrm{cofib}(\alpha) \to \mathbb{Z}[1]$, which is subsequently induced by the top-right part of the previous Cartesian square. Since the map $\beta$ is an equivalence, it follows that the map $\gamma$ is an equivalence as well.

Construction 3.13 allows us to simplify the map $\partial_R$ in Construction 3.3.

**Theorem 3.15.** *Let $R$ be a derived ring. Then the map $\partial_R$ in Construction 3.3 can be identified with the composite map*

$$L_{R/\mathbb{Z}} \otimes^{\mathbb{L}}_{R,\psi_R} \overline{R} \longrightarrow L_{\overline{R}/\mathbb{Z}} \longrightarrow L_{\overline{R}/R} \longrightarrow \overline{R}[1] \tag{3.2}$$

*in $D(\overline{R})$, where the first map is part of the transitivity sequence of $\mathbb{Z} \to R \xrightarrow{\psi_R} \overline{R}$, and the last map is classified by the map*

$$\overline{R} = R \otimes^{\mathbb{L}}_{\mathbb{Z}} \mathbb{F}_p = R \otimes^{\mathbb{L}}_{\mathbb{Z}} \mathbb{Z} \otimes^{\mathbb{L}}_{\mathbb{Z}} \mathbb{F}_p \xrightarrow{\mathrm{id}_R \otimes 1 \otimes \mathrm{id}_{\mathbb{F}_p}} R \otimes^{\mathbb{L}}_{\mathbb{Z}} \mathbb{F}_p \otimes^{\mathbb{L}}_{\mathbb{Z}} \mathbb{F}_p = \overline{R} \otimes^{\mathbb{L}}_{\mathbb{Z}} \mathbb{F}_p = \overline{R} \oplus \overline{R}[1]$$

*of derived rings. Moreover, the last map $L_{\overline{R}/R} \to \overline{R}[1]$ in (3.2) is an equivalence.*

**Proof.** First, it follows from Remark 3.14 that the last map $L_{\overline{R}/R} \to \overline{R}[1]$ in (3.2) is an equivalence. Let $R$ be a derived ring. Now we analyze the composite map

$$R \xrightarrow{\mathrm{id}_R \otimes 1 \otimes 1} R \otimes^{\mathbb{L}}_{\mathbb{Z}} \mathbb{F}_p \otimes^{\mathbb{L}}_{\mathbb{Z}} \mathbb{F}_p \xrightarrow{\varphi_{R \otimes^{\mathbb{L}}_{\mathbb{Z}} \mathbb{F}_p \otimes^{\mathbb{L}}_{\mathbb{Z}} \mathbb{F}_p}} R \otimes^{\mathbb{L}}_{\mathbb{Z}} \mathbb{F}_p \otimes^{\mathbb{L}}_{\mathbb{Z}} \mathbb{F}_p \tag{3.3}$$

classified by the map $\partial_R$ in Construction 3.3. The derived $\mathbb{F}_p$-algebra $R \otimes^{\mathbb{L}}_{\mathbb{Z}} \mathbb{F}_p \otimes^{\mathbb{L}}_{\mathbb{Z}} \mathbb{F}_p$ can be rewritten as $\overline{R} \otimes^{\mathbb{L}}_{\mathbb{F}_p, 1 \otimes \mathrm{id}} (\mathbb{F}_p \otimes^{\mathbb{L}}_{\mathbb{Z}} \mathbb{F}_p)$, where the $\mathbb{F}_p$-algebra structure on $\mathbb{F}_p \otimes^{\mathbb{L}}_{\mathbb{Z}} \mathbb{F}_p$ is given by the second factor, and the tensor product is also over it. Then the Frobenius map $\varphi_{R \otimes^{\mathbb{L}}_{\mathbb{Z}} \mathbb{F}_p \otimes^{\mathbb{L}}_{\mathbb{Z}} \mathbb{F}_p}$ is identified with the map $\varphi_{\overline{R}} \otimes \varphi_{\mathbb{F}_p \otimes^{\mathbb{L}}_{\mathbb{Z}} \mathbb{F}_p}$. Now applying Construction 3.13, we see that the composite map (3.3) can be rewritten as

$$R \xrightarrow{\psi_R} R \otimes^{\mathbb{L}}_{\mathbb{Z}} \mathbb{F}_p = R \otimes^{\mathbb{L}}_{\mathbb{Z}} \mathbb{Z} \otimes^{\mathbb{L}}_{\mathbb{Z}} \mathbb{F}_p \xrightarrow{\mathrm{id}_R \otimes 1 \otimes \mathrm{id}_{\mathbb{F}_p}} R \otimes^{\mathbb{L}}_{\mathbb{Z}} \mathbb{F}_p \otimes^{\mathbb{L}}_{\mathbb{Z}} \mathbb{F}_p.$$

The result then follows. □

**Remark 3.16.** In retrospect, it is conceptually better to interpret the composite of the first two maps in (3.2) as the map induced by the commutative diagram

$$\begin{array}{ccc} \mathbb{Z} & \longrightarrow & R \\ \downarrow & & \downarrow \\ R & \xrightarrow{\psi_R} & \overline{R} \end{array}$$

of derived rings, where the right vertical map is the canonical map $\mathrm{can}: R \to \overline{R}$, via Construction 3.22, as in the relative case (Remark 3.23).

**Remark 3.17.** We can reproduce Corollary 3.8 via Theorem 3.15 as follows. Let $R$ be a derived $\delta$-ring. Then the first map in (3.2) coincides with the composite map

$$L_{R/\mathbb{Z}} \otimes^{\mathbb{L}}_{R,\psi_R} \overline{R} \longrightarrow L_{R/\mathbb{Z}} \otimes^{\mathbb{L}}_{R,\mathrm{can}} \overline{R} \longrightarrow L_{\overline{R}/\mathbb{Z}}$$

in $D(\overline{R})$, where the first map is induced by $\varphi_R$. Now the result follows from the null-homotopy of the transitivity sequence

$$L_{R/\mathbb{Z}} \otimes^{\mathbb{L}}_{R,\mathrm{can}} \overline{R} \longrightarrow L_{\overline{R}/\mathbb{Z}} \longrightarrow L_{\overline{R}/R}$$

associated to the maps $\mathbb{Z} \to R \xrightarrow{\mathrm{can}} \overline{R}$ of derived rings.

Inspired by this description of Frobenius–Witt cotangent complex and prismatic cohomology relative to $\delta$-rings introduced in [AKN23], we introduce the following relative version of $\partial_R$ in Construction 3.3.

**Construction 3.18.** Let $A$ be a derived $\delta$-ring, and $R$ a derived $A$-algebra. Denote $\overline{R} := R \otimes_{\mathbb{Z}}^{\mathbb{L}} \mathbb{F}_p$ Then the map $\psi_R : R \to \overline{R}$ of derived rings fits into a commutative diagram

$$\begin{array}{ccc} A & \xrightarrow{\varphi_A} & A \\ \downarrow & & \downarrow \\ R & \xrightarrow{\psi_R} & \overline{R} \end{array}$$

of derived rings, and thus the composite map

$$R \xrightarrow{\text{can}} \overline{R} \xrightarrow{\psi_{\overline{R}}} \overline{R} \otimes_{\mathbb{Z}}^{\mathbb{L}} \mathbb{F}_p = \overline{R} \oplus \overline{R}[1]$$

of derived rings (here we are invoking Construction 2.2) carries a canonical $A$-algebra structure (where the $A$-algebra structure on the target is forgotten along $\varphi_A$). Thus it induces a map

$$\partial_{R/A} : L_{R/A} \otimes_{R,\psi_R}^{\mathbb{L}} \overline{R} \longrightarrow \overline{R}[1]$$

in $D(\overline{R})$. We also have a map $m_{R/A} : \overline{R} \to \operatorname{fib}(\partial_{R/A})$ in $D(\overline{R})$, along with a fiber sequence

$$\overline{R} \xrightarrow{m_{R/A}} \operatorname{fib}(\partial_{R/A}) \longrightarrow L_{R/A} \otimes_{R,\psi_R}^{\mathbb{L}} \overline{R}.$$

This leads to Frobenius–Witt cotangent complex relative to derived $\delta$-rings.

**Definition 3.19.** *Let $A$ be a derived $\delta$-ring, and $R$ a derived $A$-algebra with $\overline{R} := R \otimes_{\mathbb{Z}}^{\mathbb{L}} \mathbb{F}_p$.*

- *The* Frobenius–Witt cotangent complex $FL_{R/A}$ relative to $A$ *is the pointed $\overline{R}$-module $(\operatorname{fib}(\partial_{R/A}), m_{R/A})$ as in Construction 3.18.*
- *Let $M$ be an $\overline{R}$-module, and $m : \overline{R} \to M$ a map of $\overline{R}$-modules. Then the anima* $\operatorname{FWDer}_A(R, (M, m))$ *of* (pointed) Frobenius–Witt derivations relative to $A$ *is defined to be*

$$\operatorname{Hom}_{\mathrm{DAlg}_{A/\cdot/R}}(R, W_2(R, M, m)),$$

*where the map $W_2(R, M, m) \to R$ is the canonical map* res *in Definition 2.9, and the $A$-algebra structure on $W_2(R, M, m)$ is induced by the $A$-algebra structure on $R$, and forgetting the $A$-algebra structure on $\overline{R} \oplus \operatorname{cofib}(m) \to \overline{R} \oplus \overline{R}[1]$ along $\varphi_A : A \to A$.*

When everything is concentrated in degree 0, we get a down-to-earth description of Frobenius–Witt derivations relative to $\delta$-rings.

**Example 3.20.** Let $A$ be a $\delta$-ring, and $R$ a commutative $A$-algebra, and $M$ a $R/p$-module, and $m \in M$ an element. We see that the map $A \to W_2(R, M, m)$ of commutative rings is concretely given by $a \mapsto (a, \delta(a)\, m)$. Thus a Frobenius–Witt derivation $R \to M$ relative to $A$, which maps $p$ to $m$, can be identified with a Frobenius–Witt derivation $d : R \to M$, which maps $p$ to $m$, such that $d(a) = \delta(a)\, d(p)$ for every $a \in A$.

Applying the proof of Theorem 3.6, we see that

**Proposition 3.21.** *Let $A$ be a derived $\delta$-ring, $R$ be a derived $A$-algebra with $\overline{R} := R \otimes_{\mathbb{Z}}^{\mathbb{L}} \mathbb{F}_p$. Then the functor*

$$\begin{array}{rcl} D(\overline{R})_{\overline{R}/} & \longrightarrow & \mathrm{An} \\ (M, m) & \longmapsto & \operatorname{FWDer}_A(R, (M, m)) \end{array}$$

*is co-represented by* $(\operatorname{fib}(\partial_{R/A}), m_{R/A})$ *in Construction 3.18.*

**Construction 3.22.** Let

$$\begin{array}{ccc} A & \longrightarrow & B \\ \downarrow & & \downarrow \\ R & \longrightarrow & S \end{array} \tag{3.4}$$

be a commutative diagram of derived rings. Then we get a map

$$c: L_{R/A} \otimes^{\mathbb{L}}_{R} S \longrightarrow L_{S/B}$$

of $S$-modules. When the commutative square (3.4) is coCartesian, the map $c$ is an equivalence.

**Remark 3.23.** Let $A$ be a derived $\delta$-ring, and $R$ a derived $A$-algebra with $\overline{R} := R \otimes^{\mathbb{L}}_{\mathbb{Z}} \mathbb{F}_p$. Then as in Theorem 3.15, the map $\partial_{R/A}$ can be identified with the composite map

$$L_{R/A} \otimes^{\mathbb{L}}_{R,\psi_R} \overline{R} \longrightarrow L_{\overline{R}/A} \longrightarrow L_{\overline{R}/R} \xrightarrow{\simeq} \overline{R}[1],$$

where the first two maps are induced as in Construction 3.22 by the two commutative squares

$$\begin{array}{ccccc} A & \xrightarrow{\varphi_A} & A & \longrightarrow & R \\ \downarrow & & \downarrow & & \downarrow \\ R & \xrightarrow{\psi_R} & \overline{R} & \xrightarrow{\mathrm{id}_{\overline{R}}} & \overline{R} \end{array} \tag{3.5}$$

of derived rings, where the middle and the rightmost vertical maps are canonical maps (no Frobenius involved). In short, the map $\partial_{R/A}$ can be identified with the map in Construction 3.22 associated to the outer square of (3.5).

We give a relative and derived generalization of transitivity sequence of Frobenius–Witt cotangent complex in [Shi26, Thm 5.3].

**Construction 3.24.** Let $A$ be a derived $\delta$-ring, and $R \to S$ a map of derived $A$-algebras, with $\overline{R} := R \otimes^{\mathbb{L}}_{\mathbb{Z}} \mathbb{F}_p$ and $\overline{S} := S \otimes^{\mathbb{L}}_{\mathbb{Z}} \mathbb{F}_p$. Then by Remark 3.23 we get a commutative diagram

$$\begin{array}{ccc} L_{R/A} \otimes^{\mathbb{L}}_{R,\psi_R} \overline{R} \otimes^{\mathbb{L}}_{\overline{R}} \overline{S} & \xrightarrow{\partial_{R/A} \otimes^{\mathbb{L}}_{\overline{R}} \overline{S}} & L_{\overline{R}/R} \otimes^{\mathbb{L}}_{\overline{R}} \overline{S} \\ \downarrow & & \downarrow \simeq \\ L_{S/A} \otimes^{\mathbb{L}}_{S,\psi_S} \overline{S} & \xrightarrow{\partial_{S/A}} & L_{\overline{S}/S} \end{array}$$

in $D(\overline{S})$. The "octahedral" axiom gives rise to a fiber sequence

$$FL_{R/A} \otimes^{\mathbb{L}}_{\overline{R}} \overline{S} \longrightarrow FL_{S/A} \longrightarrow L_{S/R} \otimes^{\mathbb{L}}_{S,\psi_S} \overline{S}$$

in $D(\overline{S})$, called the *transitivity sequence* of (relative) Frobenius–Witt cotangent complex.

**Lemma 3.25.** *Let $A \to B$ be a map of animated $\delta$-rings such that the induced map $A \otimes^{\mathbb{L}}_{\mathbb{Z}} \mathbb{F}_p \to B \otimes^{\mathbb{L}}_{\mathbb{Z}} \mathbb{F}_p$ of animated $\mathbb{F}_p$-algebras is relatively perfect. Let $R$ be a derived $B$-algebra with $\overline{R} := R \otimes^{\mathbb{L}}_{\mathbb{Z}} \mathbb{F}_p$. Then the canonical map $FL_{R/A} \to FL_{R/B}$ in $D(\overline{R})_{\overline{R}/}$ is an equivalence.*

*In particular, let $A$ be a perfect $\delta$-ring, and $R$ a derived $A$-algebra. Then the Frobenius–Witt cotangent complex $FL_{R/A}$ relative to $A$ coincides with the absolute Frobenius–Witt cotangent complex $FL_R$.*

**Proof.** Since the map $A \otimes^{\mathbb{L}}_{\mathbb{Z}} \mathbb{F}_p \to B \otimes^{\mathbb{L}}_{\mathbb{Z}} \mathbb{F}_p$ of animated $\mathbb{F}_p$-algebras is relatively perfect, the cotangent complex $L_{B/A}$ is $p$-completely zero ([GR03, Lem 6.5.13(i)], cf. [Mao25, Lem 3.21] or the proof of [BF26, Lem 2.19]), and thus the canonical map $L_{R/A} \to L_{R/B}$ in $D(R)$ is $p$-completely an equivalence. The result then follows from Remark 3.23, as $\overline{R}$ is a derived $\mathbb{F}_p$-algebra. □

[Shi26, Thm 6.8] establishes vanishing of Frobenius–Witt cotangent complex of perfectoid rings. Thanks to the concrete description in Remark 3.23, we have a relative generalization (it is a generalization by Lemma 3.25).

**Theorem 3.26.** *Let $(A, I)$ be a derived prism, and $R := A/I$ with $\overline{R} := R \otimes^{\mathbb{L}}_{\mathbb{Z}} \mathbb{F}_p$. Then the map*

$$\partial_{R/A}: L_{R/A} \otimes^{\mathbb{L}}_{R,\psi_R} \overline{R} \longrightarrow \overline{R}[1]$$

*of $\overline{R}$-modules is an equivalence. Consequently, the Frobenius–Witt cotangent complex $FL_{R/A}$ relative to $A$ vanishes.*

**Proof.** By Remark 3.23, it suffices to see that the composite commutative square

$$\begin{array}{ccccc} A & \xrightarrow{\varphi_A} & A & \longrightarrow & R \\ \downarrow & & \downarrow & & \downarrow \\ R & \xrightarrow{\psi_R} & \overline{R} & \xrightarrow{\mathrm{id}} & \overline{R} \end{array}$$

in DAlg is coCartesian. This is precisely distinguishedness of the derived pre-prism $I \to A$: the canonical map $A/^{\mathbb{L}}(I, \varphi^*(I)) \to A/^{\mathbb{L}}(p, I)$ is an equivalence (cf. [Mao24, §2.4], a derived generalization of [BL22, Rem 2.1.5]). The argument uses only distinguishedness, and therefore applies more generally to distinguished derived pre-prisms in the sense of [Mao24, §2.4]. □

# 4 Logarithmic Frobenius–Witt cotangent complex

For a pre-log ring $(R, P_R, (\alpha_R: P_R \to (R, \cdot)))$ and an $R/p$-module $M$, a *logarithmic Frobenius–Witt derivation* of $(R, P_R, \alpha_R)$ valued in $M$, defined in [Tak26, Def 3.2], is a pair $(d, d^\flat)$ of maps, where $d: R \to M$ is a Frobenius–Witt derivation, and $d^\flat: P_R \to M$ is a map of abelian monoids satisfying

$$d(\alpha_R(a)) = \alpha_R(a)^p \, d^\flat(a)$$

for every $a \in P_R$. Recall that an *animated pre-log ring* (cf. [BLPØ23, §2.6]) is a triple $(R, P_R, \alpha_R)$ where $R$ is an animated ring, $P_R$ is an animated abelian monoid, and $\alpha_R: P_R \to (\Omega^\infty R, \cdot)$ a map of animated abelian monoids. We denote by $\mathrm{LAlg}^{\mathrm{an}}$ the $\infty$-category of animated pre-log rings.

We explain that previous considerations lead to a "concrete" definition of logarithmic Frobenius–Witt cotangent complex of animated pre-log rings as well. First, the arithmetic extension is canonically equipped with a pre-log structure.

**Definition 4.1.** *Let $(R, P_R, \alpha_R)$ be a pre-log ring, and $M$ an $R/p$-module, and $m \in M$ an element. The* arithmetic extension *of $(R, P_R, \alpha_R)$ by $M$ is the pre-log ring $(W_2(R, M, m), P_R \oplus M, \beta_R)$, where the map $\beta_R: P_R \oplus M \to (W_2(R, M, m), \cdot)$ induced by the map*

$$\begin{array}{ccc} P_R & \longrightarrow & W_2(R, M, m) \\ a & \longmapsto & (\alpha_R(a), 0) \end{array}$$

*and the map*

$$\begin{array}{ccc} M & \longrightarrow & W_2(R, M, m) \\ x & \longmapsto & (1, x) \end{array}$$

*of abelian monoids.*

**Remark 4.2.** Let $(R, P_R, \alpha_R)$ be a pre-log ring, and $M$ an $R/p$-module, and $m \in M$ an element. Then as in [Shi26, Prop 2.8], the map

$$\begin{array}{ccc} \left\{ \begin{array}{c} \text{pre-log ring maps } s: R \to W_2(R, M, m) \\ \text{with } \mathrm{res} \circ s = \mathrm{id} \end{array} \right\} & \to & \left\{ \begin{array}{c} \text{log FW-derivations } (d, d^\flat): R \to M \\ \text{with } d(p) = m \end{array} \right\} \\ (s_R, s_P) & \longmapsto & (\mathrm{pr}_M \circ s_R, \mathrm{pr}_M \circ s_P) \end{array}$$

is a bijection, where the "log" part of maps of pre-log rings is, by abuse of notation, omitted.

There is an extension of Example 2.1.

**Example 4.3.** Let $(R, P_R, \alpha_R)$ be a pre-log $\mathbb{F}_p$-algebra. Then the arithmetic extension $W_2(R, R, 1)$ is given by the pre-log ring $(W_2(R), P_R \oplus R, \beta_R)$ where $\beta_R$ is the map $P_R \oplus R \to W_2(R)$, $(a, x) \mapsto (\alpha_R(a), \alpha_R(a)^p x)$. By definition, this pre-log ring receives a map from the pre-log ring $(W_2(R), P_R, [\alpha_R(-)])$.

We have a pullback description of $(W_2(R, M, m), P_R \oplus M, \beta_R)$ as in Remark 2.10. For this, we have to construct the analogue $\psi_{(R, P_R, \alpha_R)}$ for animated pre-log $\mathbb{F}_p$-algebras $(R, P_R, \alpha_R)$.

**Construction 4.4.** Let $(R, P_R, \alpha_R)$ be an animated pre-log $\mathbb{F}_p$-algebra. Then the animated ring $R \otimes_{\mathbb{Z}}^{\mathbb{L}} \mathbb{F}_p$ is equivalent to the trivial square-zero extension $R \oplus R[1]$ (**not** as an animated $\mathbb{F}_p$-algebra), thus it is the trivial square zero extension $(R \oplus R[1], P_R \oplus \Omega^\infty(R[1]))$[4.1] of the pre-log ring $(R, P_R, \alpha_R)$ by the $R$-module $R[1]$, which can be viewed as an animated pre-log $\mathbb{F}_p$-algebra (**no longer** a trivial square-zero extension as such), denoted by $(R \otimes_{\mathbb{Z}}^{\mathbb{L}} \mathbb{F}_p, P_R \oplus \Omega^\infty(R[1]))$.

**Construction 4.5.** Let $(R, P_R, \alpha_R)$ be an animated pre-log $\mathbb{F}_p$-algebra. Then the map $\psi_{(R, P_R, \alpha_R)}$ of animated pre-log rings is defined to be the composite map

$$(R, P_R, \alpha_R) \xrightarrow{\text{can}} (R \otimes_{\mathbb{Z}}^{\mathbb{L}} \mathbb{F}_p, P_R \oplus \Omega^\infty(R[1])) \xrightarrow{\varphi_{(R \otimes_{\mathbb{Z}}^{\mathbb{L}} \mathbb{F}_p, P_R \oplus \Omega^\infty(R[1]))}} (R \otimes_{\mathbb{Z}}^{\mathbb{L}} \mathbb{F}_p, P_R \oplus \Omega^\infty(R[1]))$$

of animated pre-log rings, where the Frobenius map on $P_R \oplus \Omega^\infty(R[1])$ is multiplication by $p$. Note that the multiplication by $p$ on $\Omega^\infty(R[1])$ is null-homotopic, and the null-homotopy is functorially given by the map $\mathbb{F}_p \to R$ of animated rings.

**Proposition 4.6.** *Let $(R, P_R, \alpha_R)$ be a pre-log ring with $\overline{R} := R \otimes_{\mathbb{Z}}^{\mathbb{L}} \mathbb{F}_p$, and $M$ an $R/p$-module, and $m \in M$ an element. Then the arithmetic extension $(W_2(R, M, m), P_R \oplus M, \beta)$ of $(R, P_R, \alpha_R)$ by $M$ coincides with the pullback $W_2^{\text{pb}}$ of the diagram*

$$\begin{array}{ccc} & & (R, P_R, \alpha_R) \\ & & \downarrow \text{can} \\ & & (\overline{R}, P_R, \text{can} \circ \alpha_R) \\ & & \downarrow \psi_{(\overline{R}, P_R)} \\ (\overline{R} \oplus \text{cofib}(m), P_R \oplus \Omega^\infty \text{cofib}(m)) & \longrightarrow & (\overline{R} \otimes_{\mathbb{Z}}^{\mathbb{L}} \mathbb{F}_p, P_R \oplus \Omega^\infty(\overline{R}[1])) \end{array}$$

*of pre-log rings, where the bottom left is the trivial square-zero extension of $(\overline{R}, P_R, \text{can} \circ \alpha_R)$ by* $\text{cofib}(m)$*, where $m$ is viewed as a map $\overline{R} \to M$ of $\overline{R}$-modules*[4.2].

**Proof.** By Remark 2.11, the pullback $W_2^{\text{pb}}$ has underlying ring $W_2(R, M, m)$. Moreover, the underlying abelian monoid is the pullback of the diagram

$$\begin{array}{ccc} & & P_R \\ & & \downarrow \\ P_R \oplus \Omega^\infty \text{cofib}(m) & \longrightarrow & P_R \oplus \Omega^\infty(\overline{R}[1]) \end{array}$$

of animated abelian monoids, where the vertical map is the composite map

$$P_R \xrightarrow{p} P_R \xhookrightarrow{\text{id}_{P_R}} P_R \oplus \Omega^\infty(\overline{R}[1]),$$

see Construction 4.5. It follows that the above pullback is the animated monoid $P_R \oplus \Omega^\infty M$, which is concentrated in degree 0. This shows that the pullback $W_2^{\text{pb}}$ is a classical pre-log ring.

It remains to check that the exponential map $e^{\text{pb}}: P_R \oplus M \to (W_2(R, M, m), \cdot)$ in $W_2^{\text{pb}}$ coincides with the map $\beta_R: P_R \oplus M \to (W_2(R, M, m), \cdot)$. Tracing through the identification in Remark 2.11 and Proposition 2.6, one can see that the map $e^{\text{pb}}|_{P_R}$ is the composite

$$P_R \xrightarrow{\alpha_R} (R, \cdot) \xrightarrow{[-]} (W_2(R, M, m), \cdot)$$

4.1. We omit the map $P_R \oplus \Omega^\infty(R[1]) \to (\Omega^\infty(R \oplus R[1]), \cdot)$ in the notation.

4.2. Since $M$ is a $\pi_0(\overline{R})$-module, it is also an $\overline{R}$-module.

of abelian monoids[4.3], which is the same as $\beta_R|_{P_R}$. Now we identify $e^{\mathrm{pb}}|_M$, which by definition, appears in the map

$$\begin{array}{ccccc} M & \xrightarrow{e^{\mathrm{pb}}|_M} & \mathrm{GL}_1(W_2(R,M,m)) & \xrightarrow{\mathrm{res}} & \mathrm{GL}_1(R) \\ \| & & \downarrow & & \downarrow \\ M & \xrightarrow{1+(-)} & \mathrm{GL}_1(\overline{R}\oplus \mathrm{cofib}(m)) & \longrightarrow & \mathrm{GL}_1(\overline{R}\otimes^{\mathbb{L}}_{\mathbb{Z}}\mathbb{F}_p) \end{array}$$

of fiber sequences in $D(\mathbb{Z})$ (and everything in question is actually connective). Comparing this with the map

$$\begin{array}{ccccc} M & \xrightarrow{V} & W_2(R,M,m) & \xrightarrow{\mathrm{res}} & R \\ \| & & \downarrow & & \downarrow \\ M & \longrightarrow & \overline{R}\oplus \mathrm{cofib}(m) & \longrightarrow & \overline{R}\otimes^{\mathbb{L}}_{\mathbb{Z}}\mathbb{F}_p \end{array}$$

of fiber sequences in $D(\mathbb{Z})$ (see Remark 2.13), we see that $e^{\mathrm{pb}}|_M = 1 + V = \beta_R|_M$. □

Thus we extend the arithmetic extension and logarithmic Frobenius–Witt derivations to general animated pre-log rings.

**Definition 4.7.** *Let $(R,P_R,\alpha_R)$ be an animated pre-log ring with $\overline{R} := R\otimes^{\mathbb{L}}_{\mathbb{Z}}\mathbb{F}_p$, and $M$ a connective $\overline{R}$-module, and $m:\overline{R}\to M$ a map of $\overline{R}$-modules. Then we define*

- *the* arithmetic extension *of animated pre-log ring of $(R,P_R,\alpha_R)$ by $(M,m)$ to be the pullback of the diagram*

$$\begin{array}{ccc} & & (R,P_R,\alpha_R) \\ & & \downarrow{\scriptstyle \mathrm{can}} \\ & & (\overline{R},P_R,\mathrm{can}\circ\alpha_R) \\ & & \downarrow{\scriptstyle \psi_{(\overline{R},P_R)}} \\ (\overline{R}\oplus\mathrm{cofib}(m), P_R\oplus\Omega^\infty\mathrm{cofib}(m)) & \longrightarrow & (\overline{R}\otimes^{\mathbb{L}}_{\mathbb{Z}}\mathbb{F}_p, P_R\oplus\Omega^\infty(\overline{R}[1])) \end{array}$$

  *of animated pre-log rings.*

- *the anima* $\mathrm{LFWDer}((R,P_R,\alpha_R),(M,m))$ *of* logarithmic Frobenius–Witt derivations *is defined to be*

$$\mathrm{Hom}_{\mathrm{LAlg}^{\mathrm{an}}_{/(R,P_R,\alpha_R)}}((R,P_R,\alpha_R),(W_2(R,M,m),P_R\oplus\Omega^\infty M,\beta_R)).$$

Arguing as in Theorem 3.6, we see that

**Theorem 4.8.** *Let $(R,P_R,\alpha_R)$ be an animated pre-log ring with $\overline{R} := R\otimes^{\mathbb{L}}_{\mathbb{Z}}\mathbb{F}_p$, and $M$ a connective $\overline{R}$-module, and $m:\overline{R}\to M$ a map of $\overline{R}$-modules. Then the functor*

$$\begin{array}{rcl} D_{\geqslant 0}(\overline{R})_{\overline{R}/} & \longrightarrow & \mathrm{An} \\ (M,m) & \longmapsto & \mathrm{LFWDer}((R,P_R,\alpha_R),(M,m)) \end{array}$$

*is co-represented by $(\mathrm{fib}(\partial_{(R,P_R,\alpha_R)}), m_{(R,P_R,\alpha_R)})$, where $m_{(R,P_R,\alpha_R)}$ is the fiber of the map*

$$\partial_{(R,P_R,\alpha_R)} : L_{(R,P_R,\alpha_R)/\mathbb{Z}}\otimes^{\mathbb{L}}_{R,\psi_R}\overline{R}\longrightarrow \overline{R}[1]$$

4.3. It is essentially deduced from the fact that, for every animated ring $R$ (or equivalently, every finitely generated polynomial ring $R$), the Teichmüller map $[-]:(\Omega^\infty R,\cdot)\to(\Omega^\infty W_2(R),\cdot)$ is given by the diagram

$$\begin{array}{ccc} (\Omega^\infty R,\cdot) & \xrightarrow{\mathrm{id}} & (\Omega^\infty R,\cdot) \\ \downarrow{\scriptstyle p} & & \downarrow{\scriptstyle \Omega^\infty\psi_R} \\ (\Omega^\infty R,\cdot) & \xrightarrow{\Omega^\infty(\mathrm{can})} & (\Omega^\infty(R\otimes^{\mathbb{L}}_{\mathbb{Z}}\mathbb{F}_p),\cdot) \end{array}$$

of animated abelian monoids, where we realize $W_2(R)$ as the pullback $R\times_{\mathrm{can},R\otimes^{\mathbb{L}}_{\mathbb{Z}}\mathbb{F}_p,\psi_R}R$.

*of $\overline{R}$-modules classifying the composite map*

$$(R, P_R, \alpha_R) \xrightarrow{\text{can}} (\overline{R}, P_R, \text{can} \circ \alpha_R) \xrightarrow{\psi_{\overline{R}}} (\overline{R} \oplus \overline{R}[1], P_R \oplus \Omega^\infty(\overline{R}[1]))$$

*of animated pre-log rings into the trivial square-zero extension of* $(\overline{R}, P_R, \text{can} \circ \alpha_R)$ *by* $\overline{R}[1]$.

Similarly to Theorem 3.15 and Remark 3.16, we have

**Theorem 4.9.** *Let* $(R, P_R, \alpha_R)$ *be an animated pre-log ring with* $\overline{R} := R \otimes_{\mathbb{Z}}^{\mathbb{L}} \mathbb{F}_p$. *Then the map* $\partial_{(R,P_R,\alpha_R)}$ *can be identified with the map*

$$L_{(R,P_R,\alpha_R)/\mathbb{Z}} \otimes_{R,\psi_R}^{\mathbb{L}} \overline{R} \longrightarrow L_{(\overline{R},P_R,\text{can}\circ\alpha_R)/(R,P_R,\alpha_R)}$$

*induced by the commutative diagram*

$$\begin{array}{ccc} \mathbb{Z} & \longrightarrow & (R, P_R, \alpha_R) \\ \downarrow & & \downarrow \\ (R, P_R, \alpha_R) & \xrightarrow{\psi_R} & (\overline{R}, P_R, \text{can} \circ \alpha_R) \end{array}$$

*of animated pre-log rings, as in Construction 3.22.*

# 5 A regularity criterion

In this section, independently of previous developments, we give a regularity criterion of Noetherian rings without $F$-finiteness or locality (Theorem 5.2). This is inspired by [BM23, Lem 4.18].

**Lemma 5.1.** *Let* $(R, \mathfrak{m}, k)$ *be a local ring with residue field* $k$ *of char* $p$, *with* $\overline{R} := R \otimes_{\mathbb{Z}}^{\mathbb{L}} \mathbb{F}_p$.

- *Suppose that the Frobenius–Witt cotangent complex* $FL_R$ *is flat as an* $\overline{R}$*-module. Then the cotangent complex* $L_{k/R}$ *is concentrated in degrees* $\leq 1$.
- *Conversely, suppose that the cotangent complex* $L_{k/R}$ *is concentrated in degrees* $\leq 1$. *Then the* $k$*-module* $FL_R \otimes_{\overline{R}}^{\mathbb{L}} k$ *is concentrated in degrees* $\leq 0$, *where the map* $\overline{R} \to k$ *of animated* $\mathbb{F}_p$*-algebras is induced by the canonical map* $R \to k$ *and the* $\mathbb{F}_p$*-algebra structure on* $k$.

**Proof.** Let $\bar{k} := k \otimes_{\mathbb{Z}}^{\mathbb{L}} \mathbb{F}_p$. We examine the transitivity sequence ([Shi26, Thm 5.3] or Construction 3.24)

$$FL_R \otimes_{\overline{R}}^{\mathbb{L}} \bar{k} \longrightarrow FL_k \longrightarrow L_{k/R} \otimes_{k,\psi_k}^{\mathbb{L}} \bar{k}$$

in $D(\bar{k})$ associated to $R \to k$. By [Shi26, Cor 6.4], the $\bar{k}$-module $FL_k$ is equivalent to the $\bar{k}$-module $L_{k/\mathbb{F}_p} \otimes_{k,\psi_k}^{\mathbb{L}} \bar{k}$. We base change this fiber sequence along the map $\bar{k} = k \otimes_{\mathbb{Z}}^{\mathbb{L}} \mathbb{F}_p \to k$ induced by $\mathrm{id}_k$ and the $\mathbb{F}_p$-algebra structure on $k$, we get a fiber sequence

$$FL_R \otimes_{\overline{R}}^{\mathbb{L}} k \longrightarrow L_{k/\mathbb{F}_p} \otimes_{k,\varphi_k}^{\mathbb{L}} k \longrightarrow L_{k/R} \otimes_{k,\varphi_k}^{\mathbb{L}} k \tag{5.1}$$

in $D(k)$. Since the $k$-module $L_{k/\mathbb{F}_p}$ is actually flat [BLM18, Prop 9.5.1], so is the $k$-module $L_{k/\mathbb{F}_p} \otimes_{k,\varphi_k}^{\mathbb{L}} k$.

Now we suppose that the $\overline{R}$-module $FL_R$ is flat. Then the fiber sequence (5.1) implies that the $k$-module $L_{k/R} \otimes_{k,\varphi_k}^{\mathbb{L}} k$ has Tor-amplitude $\leq 1$. Since $k$ is a field, the Frobenius map $\varphi_k : k \to k$ is faithfully flat. It follows that the $k$-module $L_{k/R}$ has Tor-amplitude $\leq 1$, or equivalently, is concentrated in degrees $\leq 1$.

Conversely, suppose that the cotangent complex $L_{k/R}$ is concentrated in degrees $\leq 1$, then so is the $k$-module $L_{k/R} \otimes_{k,\varphi_k}^{\mathbb{L}} k$ by flatness of $\varphi_k : k \to k$. It follows from the fiber sequence (5.1) that the $k$-module $FL_R \otimes_{\overline{R}}^{\mathbb{L}} k$ is concentrated in degrees $\leq 0$. □

**Theorem 5.2.** *Let $R$ be a Noetherian ring with $\overline{R} := R \otimes_{\mathbb{Z}}^{\mathbb{L}} \mathbb{F}_p$.*

- *Suppose that the ring $R$ is $p$-local (i.e. $p \in \operatorname{Rad}(R)$), and the $\overline{R}$-module $FL_R$ is flat. Then the Noetherian ring $R$ is regular.*
- *Conversely, suppose that the ring $R$ is regular. Then the $\overline{R}$-module $FL_R$ is flat.*

**Proof.** For every prime ideal $\mathfrak{p} \in \operatorname{Spec}(R/p)$, we will denote by $\kappa(\mathfrak{p})$ the residue field $R_{\mathfrak{p}}/\mathfrak{p} R_{\mathfrak{p}} = \operatorname{Frac}(R/\mathfrak{p})$.

Suppose that the $\overline{R}$-module $FL_R$ is flat. Then for every maximal ideal $\mathfrak{m} \subseteq R$, since $p \in \operatorname{Rad}(R) \subseteq \mathfrak{m}$ (i.e. $\mathfrak{m} \in \operatorname{Spec}(R/p)$), the residue field $\kappa(\mathfrak{m}) = R/\mathfrak{m}$ is of char $p$. Note that the $\overline{R}_{\mathfrak{m}}$-module $FL_{R_{\mathfrak{m}}}$, which is equivalent to the $\overline{R}_{\mathfrak{m}}$-module $FL_R \otimes_{\overline{R}}^{\mathbb{L}} \overline{R}_{\mathfrak{m}}$ by transitivity sequence associated to $R \to R_{\mathfrak{m}}$ (cf. [Shi26, Rem 5.5]), is flat. It follows from Lemma 5.1 that the cotangent complex $L_{\kappa(\mathfrak{m})/R_{\mathfrak{m}}}$ is concentrated in degrees $\leq 1$, and by [Iye07, Prop 8.12] (where we use Noetherianness), the Noetherian local ring $R_{\mathfrak{m}}$ is regular. It follows from [Sta21, Tag 0AFS] that the Noetherian ring $R$ is regular.

Conversely, suppose that the Noetherian ring $R$ is regular. By [BLM18, Lem 9.5.2][5.1] and connectivity of $FL_R$, it suffices to check that, for every prime ideal $\mathfrak{p} \in \operatorname{Spec}(R/p)$, the connective $\kappa(\mathfrak{p})$-module $FL_R \otimes_{\overline{R}}^{\mathbb{L}} \kappa(\mathfrak{p})$ is concentrated in degrees $\leq 0$. This follows from applying Lemma 5.1 to the regular local ring $R_{\mathfrak{p}}$. □

# Appendix A Flatness of Frobenius–Witt cotangent complex (ideas from ChatGPT)

In this section, we record some flatness results of Frobenius–Witt cotangent complex essentially due to ChatGPT-6 Astra. A crucial observation is the flatness of Frobenius–Witt cotangent complex of valuation rings (Theorem A.2), which gives a further indication that Frobenius–Witt cotangent complex might play an important role in singularity theory in mixed characteristic in addition to Theorem 5.2. Another result is an equivalent condition in terms of the cotangent complex of the Frobenius map $\psi_{(-)}$ (Corollary A.10). Before that, we record a simple lemma.

**Lemma A.1.** *Let $S$ be an animated ring whose $p$-completion $S_p^{\wedge}$ is perfectoid. Then the Frobenius–Witt cotangent complex $FL_S$ is contractible.*

**Proof.** The $p$-completion map $S \to S_p^{\wedge}$ induces an equivalence of Frobenius–Witt cotangent complexes, either by [Shi26, Cor 5.6] or our description in Theorem 3.6. Then the result follows from [Shi26, Thm 6.8] (cf. Theorem 3.26). □

Here is the main result whose proof was developed in discussions with ChatGPT-6 Astra.

**Theorem A.2.** *Let $V$ be a valuation ring, and $\overline{V} := V \otimes_{\mathbb{Z}}^{\mathbb{L}} \mathbb{F}_p$. Then the Frobenius–Witt cotangent complex $FL_V$ of $V$, as a $\overline{V}$-module, is flat.*

**Proof.** Let $K := \operatorname{Frac}(V)$ be the fractional field with valuation $|\cdot|_K$ associated to $V$, and $K^a$ an algebraic closure of $K$. We extend the valuation $|\cdot|_K$ to a valuation $|\cdot|_{K^a}$ on $K^a$, and let $V^a$ denote the associated valuation ring (which depends on the chosen extension $|\cdot|_{K^a}$, not only on $V$). Then the $p$-completion of $V^a$ is a perfectoid ring.

Denote $\overline{V^a} := V^a \otimes_{\mathbb{Z}}^{\mathbb{L}} \mathbb{F}_p$. We examine the transitivity sequence (by [Shi26, Thm 5.3], cf. Construction 3.24)

$$FL_V \otimes_{\overline{V}}^{\mathbb{L}} \overline{V^a} \longrightarrow FL_{V^a} \longrightarrow L_{V^a/V} \otimes_{V^a, \psi_{V^a}}^{\mathbb{L}} \overline{V^a}$$

in $D(\overline{V^a})$ associated to the map $V \to V^a$. By Lemma A.1, the $\overline{V^a}$-module $FL_{V^a}$ is contractible.

5.1. Since for every animated ring $R$, by definition, a connective $R$-module $M$ is flat if and only if the connective $\pi_0(R)$-module $M \otimes_R^{\mathbb{L}} \pi_0(R)$ is flat. Thus [BLM18, Lem 9.5.2] extends to animated rings.

By [GR03, Thm 6.3.32], the $V^a$-module $L_{V^a/V}$ is concentrated in degrees $[0,1]$, with the $V^a$-module $\pi_1(L_{V^a/V})$ being torsion-free. Since the $V^a$-module $L_{V^a/V}$ is an extension of $\pi_0(L_{V^a/V})$ and $\pi_1(L_{V^a/V})[1]$, both of which have Tor-amplitude $\leq 1$ as $V^a$ is a valuation ring (thus of Tor-dimension $\leq 1$, with [Sta21, Tag 0539] invoked), the $V^a$-module $L_{V^a/V}$ has Tor-amplitude $\leq 1$ as well.

Consequently, the $\overline{V^a}$-module

$$FL_V \otimes^{\mathbb{L}}_{\bar{V}} \overline{V^a} \simeq (L_{V^a/V} \otimes^{\mathbb{L}}_{V^a, \psi_{V^a}} \overline{V^a})[-1]$$

has Tor-amplitude $\leq 0$. The extension $V \to V^a$ is faithfully flat by [Sta21, Tag 0539 & Tag 00HR], which implies that the map $\bar{V} \to \overline{V^a}$ is faithfully flat as well by [Sta21, Tag 00HI], and the result then follows. □

In retrospect, the proof of Theorem A.2 is very similar to [Bou23, Prop 3.8]. Actually, Theorem A.2 implies Corollary A.5, which incorporates both the char $p$ case and the mixed characteristic case, without re-running the argument of Theorem A.2, whose proof in the mixed characteristic case was developed in discussions with ChatGPT-6 Astra. For this, we need the following definition and flatness criterion from [Hat60].

**Definition A.3.** **([Hat60, §1])** *Let $R$ be an associative ring. We say that a right $R$-module $M$ (concentrated in degree 0) is* torsion-free *if, for every $r \in R$, the canonical inclusion*

$$M \operatorname{Ker}(m_r : R \to R) \subseteq \operatorname{Ker}(m_r : M \to M)$$

*is a bijection, where $m_r$ is the right multiplication by $r$.*

**Proposition A.4.** **([Hat60, Prop 2])** *Let $R$ be an associative ring, and $M$ a right $R$-module. if the right $R$-module $M$ is flat, then it is torsion-free. The converse is true if every finitely generated left $R$-ideal is principal.*

**Corollary A.5.** *Let $V \to W$ be an extension of valuation rings. Suppose that the $p$-completion of $V$ is perfectoid. Then the $W$-module $L_{W/V}$ is $p$-completely flat.*

**Proof.** Let $\bar{V} := V \otimes^{\mathbb{L}}_{\mathbb{Z}} \mathbb{F}_p$ and $\overline{W} := W \otimes^{\mathbb{L}}_{\mathbb{Z}} \mathbb{F}_p$. We start with the transitivity sequence

$$FL_V \otimes^{\mathbb{L}}_{\bar{V}} \overline{W} \longrightarrow FL_W \longrightarrow L_{W/V} \otimes^{\mathbb{L}}_{W, \psi_W} \overline{W}$$

in $D(\overline{W})$ by [Shi26, Thm 5.3], cf. Construction 3.24. By Lemma A.1, the first term is contractible. It follows that the $\overline{W}$-module $L_{W/V} \otimes^{\mathbb{L}}_{W, \psi_W} \overline{W}$ is flat.

- When $V$ is a valuation $\mathbb{F}_p$-algebra, then we may further base-change along the multiplication map $\overline{W} = W \otimes^{\mathbb{L}}_{\mathbb{Z}} \mathbb{F}_p \to W$, deducing the flatness of $W$-module $L_{W/V} \otimes^{\mathbb{L}}_{W, \varphi_W} W$. Now as the Frobenius map $\varphi_W : W \to W$ on the valuation $\mathbb{F}_p$-algebra is faithfully flat by [Sta21, Tag 0539 & Tag 00HR], the flatness of $L_{W/V}$ then follows.
- Now we assume that $V$ is $p$-torsion-free, then so is $W$. If $p$ is invertible in $V$, then the cotangent complex is $p$-completely contractible, thus $p$-completely flat. Now we assume that $p$ is not invertible in $V$, i.e. $p$ lies in the maximal ideal $\mathfrak{m}_V \subseteq V$, and $V$ is a mixed characteristic valuation ring in [Bou23, Lem 3.5]. By Lemma A.6, we may pick $\pi \in V$ such that $\pi^p W = pW$.

  Let $W' := W/\pi$. Then the Frobenius map $\varphi_{\overline{W}}$ factors as

  $$\overline{W} \xrightarrow{\text{can}} W' \xrightarrow{\varphi'_{\overline{W}}} \overline{W} \tag{A.1}$$

  where the first map is a surjection with kernel being nilpotent (of exponent $p$). Now we show that the second map $\varphi'_{\overline{W}} : W' \to \overline{W}$ is faithfully flat. Since $W$ is a valuation ring, every finitely generated ideal of $W$ is principal, thus so is $W'$. Now we check that the $W'$-module $\overline{W}$ (via $\varphi'_{\overline{W}}$) is torsion-free. Indeed, for every $a \in W$, concerning the division relation between $a$ and $\pi$, there are two cases

  - If $a \in \pi W$, then the multiplication $m_a$ by $a$ is zero on $W'$, thus $\overline{W} \operatorname{Ker}(m_a : W' \to W') = \overline{W} = \operatorname{Ker}(m_{\varphi'_{\overline{W}}(a)} : \overline{W} \to \overline{W})$;

- Otherwise, $a \notin \pi W$ (and in particular, $a \neq 0$) and $\pi \in a W$. Then for every $x \in W$ such that $a^p x \in p W = \pi^p W$, we have $x \in (\pi / a)^p W$. This also shows that $\overline{W}\,\mathrm{Ker}(m_a : W' \to W') = \mathrm{Ker}(m_{\varphi'_{\overline{W}}(a)} : \overline{W} \to \overline{W})$.

By Proposition A.4, we see that the map $\varphi'_{\overline{W}} : W' \to \overline{W}$ of rings is flat. Now by [Sta21, Tag 0BR6], the map $\mathrm{can} : \overline{W} \to W'$ is a homeomorphism on $\mathrm{Spec}(-)$, and by [Sta21, Tag 0CCB], so is the map $\varphi_{\overline{W}}$. It follows that the map $\varphi'_{\overline{W}}$ is a homeomorphism on $\mathrm{Spec}(-)$ as well. Consequently, by [Sta21, Tag 00HQ], the map $\varphi'_{\overline{W}}$ of rings is faithfully flat.

We have already known that the $\overline{W}$-module $L_{W/V} \otimes^{\mathbb{L}}_{W, \psi_W} \overline{W}$ is flat. By [Sta21, Tag 0H75] and [Sta21, Tag 068S] along with the factorization (A.1), we see that the base change along $\varphi_{\overline{W}}$ reflects flatness of connective modules, and the result follows. □

We need the following lemma, which was summarized from discussions with ChatGPT-6 Astra.

**Lemma A.6.** *Let $S$ be a $p$-local animated ring such that its (derived) $p$-completion $S_p^\wedge$ is perfectoid. Then there exists $\pi \in \pi_0(S)$ and $u \in \pi_0(S)^\times$ such that $\pi^p = p\, u$ in $\pi_0(S)$.*

**Proof.** We pick $\hat{\pi} \in S_p^\wedge$ and $\hat{u} \in (S_p^\wedge)^\times$ such that $\hat{\pi}^p = p\, \hat{u}$ (cf. [BMS18, Rem 3.8]). Note that the $p$-completion map $S \to S_p^\wedge$ induces an equivalence

$$S \otimes^{\mathbb{L}}_{\mathbb{Z}} (\mathbb{Z} / p^2) \xrightarrow{\simeq} S_p^\wedge \otimes^{\mathbb{L}}_{\mathbb{Z}} (\mathbb{Z} / p^2)$$

of animated rings, thus there exists $\pi \in \pi_0(S)$ and $u_0 \in \pi_0(S)^\times$ whose images in $\pi_0(S \otimes^{\mathbb{L}}_{\mathbb{Z}} (\mathbb{Z} / p^2)) = \pi_0(S) / p^2$ are congruent to $\hat{\pi}$ and $\hat{u}$ respectively, thus there exists $\varepsilon \in \pi_0(S)$ such that $\pi^p - p\, u_0 = p^2 \varepsilon$. Set $u := u_0 + p\, \varepsilon$, which is invertible in $\pi_0(S)$ as $u_0 \in \pi_0(S)^\times$ and $p \in \mathrm{Rad}(\pi_0(S))$. □

**Example A.7.** In Corollary A.5, setting $V = \mathbb{F}_p$, we recover the char $p$ case [KST21, Cor A.4(i,ii)].

**Example A.8.** In Corollary A.5, let $V$ be a mixed characteristic valuation ring whose $p$-completion is perfectoid. Then by [GR03, Thm 6.5.12], the (non-$p$-completed) cotangent complex $L_{W/V}$ is concentrated in degree 0, and we recover [Bou23, Prop 3.8]: the $p$-complete flatness of $L_{W/V}$ along with the $p$-torsion-freeness of $W$ implies the $p$-torsion-freeness of $\Omega^1_{W/V}$.

The following construction relating Frobenius–Witt cotangent complex and usual cotangent complex associated to the derived Frobenius map is essentially discovered by ChatGPT-6 Astra.

**Construction A.9.** Let $R$ be an animated ring with $\overline{R} := R \otimes^{\mathbb{L}}_{\mathbb{Z}} \mathbb{F}_p$. We will realize the $\overline{R}$-module $FL_R[1]$ as a direct summand of the cotangent complex $L_{\psi_R}$ of the map $\psi_R : R \to \overline{R}$ as follows. The $\overline{R}$-module $L_{\psi_R}$ fits into the transitivity sequence

$$L_{R/\mathbb{Z}} \otimes^{\mathbb{L}}_{R, \psi_R} \overline{R} \xrightarrow{\mathrm{d}_{\psi_R}} L_{\overline{R}/\mathbb{Z}} \longrightarrow L_{\psi_R} \tag{A.2}$$

of $\overline{R}$-modules. Since the derived ring $R$ is the pushout of the animated rings $R$ and $\mathbb{F}_p$, the map

$$L_{\overline{R}/\mathbb{Z}} \longrightarrow L_{\overline{R}/R} \oplus L_{\overline{R}/\mathbb{F}_p},$$

induced by canonical maps $R \to \overline{R}$ and $\mathbb{F}_p \to \overline{R}$, is an equivalence. The composite map

$$L_{R/\mathbb{Z}} \otimes^{\mathbb{L}}_{R, \psi_R} \overline{R} \xrightarrow{\mathrm{d}_{\psi_R}} L_{\overline{R}/\mathbb{Z}} \longrightarrow L_{\overline{R}/\mathbb{F}_p}$$

in $D(\overline{R})$ can be identified with the differential

$$L_{\overline{R}/\mathbb{F}_p} \otimes^{\mathbb{L}}_{\overline{R}, \varphi_{\overline{R}}} \overline{R} \xrightarrow{\mathrm{d}_{\varphi_{\overline{R}}}} L_{\overline{R}/\mathbb{F}_p}$$

which is null-homotopic by [GR03, Lem 6.5.13(i)], cf. [Mao25, Lem 3.21] or the proof of [BF26, Lem 2.19], while the composite map

$$L_{R/\mathbb{Z}} \otimes^{\mathbb{L}}_{R, \psi_R} \overline{R} \xrightarrow{\mathrm{d}_{\psi_R}} L_{\overline{R}/\mathbb{Z}} \longrightarrow L_{\overline{R}/R}$$

is precisely the map $\partial_R$ in Construction 3.3 by Theorem 3.15. Combining the preceding observations, we may rewrite the fiber sequence (A.2) as the fiber sequence

$$L_{R/\mathbb{Z}} \otimes^{\mathbb{L}}_{R,\psi_R} \overline{R} \xrightarrow{\binom{\partial_R}{0}} L_{\overline{R}/R} \oplus L_{\overline{R}/\mathbb{F}_p} \longrightarrow L_{\psi_R}$$

in $D(\overline{R})$, which gives rise to an equivalence

$$FL_R[1] \oplus L_{\overline{R}/\mathbb{F}_p} \xrightarrow{\simeq} L_{\psi_R}$$

of $\overline{R}$-modules.

Consequently, we have

**Corollary A.10.** *Let $R$ be an animated ring with $\overline{R} := R \otimes^{\mathbb{L}}_{\mathbb{Z}} \mathbb{F}_p$. Then the following statements are equivalent.*

1. *The animated ring $R$ is $p$-quasi-lci*[A.1]*, and the $\overline{R}$-module $FL_R$ is flat.*
2. *The cotangent complex $L_{\psi_R}$ of the map $\psi_R: R \to \overline{R}$, as an $\overline{R}$-module, has* Tor*-amplitude* $\leq 1$.

**Proof.** The statement

the $\overline{R}$-module $L_{\psi_R}$ has Tor-amplitude $\leq 1$

is equivalent to the statement

both $\overline{R}$-modules $FL_R[1]$ and $L_{\overline{R}/\mathbb{F}_p}$ has Tor-amplitude $\leq 1$

by Construction A.9, which is equivalent to the statement the $\overline{R}$-module $FL_R$ is flat, and the $R$-module $L_{R/\mathbb{Z}}$ has $p$-complete Tor-amplitude $\leq 1$

by connectivity of $FL_R$ (as $R$ is connective). ☐

A.1. As for usual rings, we say that an animated ring $R$ is *$p$-quasi-lci* if the $R$-module $L_{R/\mathbb{Z}}$ has $p$-complete Tor-amplitude $\leq 1$.